\documentclass[AMA,STIX1COL]{WileyNJD-v2}
\usepackage{moreverb}
\usepackage{amsmath}
\usepackage{amsfonts}
\usepackage{amssymb}
\usepackage{amsthm}
\usepackage{graphicx}
\usepackage{array}
\usepackage{tabularx}
\usepackage{color}
\usepackage{siunitx}
\usepackage{caption}
\usepackage{subcaption}

\renewcommand{\P}{\mathbf{\hat P}}

\newcommand{\bQ}{\mathbf{Q}}
\newcommand{\bF}{\mathbf{F}}

\articletype{Article Type}%

\received{<day> <Month>, <year>}
\revised{<day> <Month>, <year>}
\accepted{<day> <Month>, <year>}

\begin{document}


\title{A flux-vector splitting scheme for the shallow water equations extended to high-order on unstructured meshes}  

\author[UNITN]{E. F. Toro}

\author[UTA]{C. E. Castro}

\author[ETH]{D. Vanzo}

\author[UNITN]{A. Siviglia*}

\authormark{Toro \textsc{et al}}

\address[UNITN]{\orgdiv{Laboratory of Applied Mathematics}, \orgname{DICAM, University of Trento}, \orgaddress{\country{Italy}}}

\address[UTA]{\orgdiv{Departamento de Ingenier{\'i}a Mec{\'a}nica, Facultad de Ingenier{\'i}a}, \orgname{Universidad de Tarapac{\'a}}, \orgaddress{\country{Chile}}}

\address[ETH]{\orgdiv{Laboratory of Hydraulics, Hydrology and Glaciology VAW}, \orgname{ETH Zurich}, \orgaddress{\country{Switzerland}}}

\corres{*Annunziato Siviglia. \email{annunziato.siviglia@unitn.it}}



\abstract[Abstract]{We present an advection-pressure flux-vector splitting method for the one and two-dimensional shallow water equations following the  approach first proposed  by Toro and V\'azquez  \cite{Toro:2012b} for the compressible Euler equations. The resulting first-order schemes turn out to be exceedingly simple, with accuracy and robustness comparable to that of the sophisticated Godunov upwind method used in conjunction with complete non-linear Riemann solvers. The technique splits the full system into two subsystems,  namely an advection system and a  pressure system. The sought numerical flux results  from fluxes for each of the subsystems. The basic methodology,  extended on 2D unstructured meshes, constitutes the building block for the construction of numerical schemes of very high order of accuracy following the ADER approach.  The presented numerical schemes are systematically assessed on a carefully selected suite of test problems with reference solutions,  in one and two space dimensions.The applicability of the schemes is illustrated through  simulations of tsunami wave propagation in the Pacific Ocean.}



\maketitle


\section{Introduction}

The rapidily evolving environmental changes resulting from the interaction of humans with surrounding natural systems is the subject of increasing concern to scientists, policy makers and the population at large.  In this scenario of accelerating geophysical changes provoking a real environmental emergency, the role of environmental sciences, in a broad sense, is of fundamental importance.   This situation calls for a multidisciplinary approach aimed, first of all,  at understanding basic mechanisms so as to be able to recommend robust measures to industry, decision makers, governments and international agencies. Mathematical modelling and computational simulation, coupled to experimental techniques, observation, data collection and measurements, plays a fundamental role here. Such role is expected to increase in the future and will include the development of efficient and reliable computational methods, based on physically meaningful mathematical models, increasing computational power and an ever increasing level of available environmental data. \\

This paper is concerned with mathematical models and numerical methods for simulating water flows in channels, rivers, lakes and oceans.  The models and methods may also be applicable to shallow water type models for atmospheric flows. Specifically, we are concerned with new numerical methods to solve the one- and two-dimensional time dependent shallow water equations. The development of numerical methods for such partial differential equations has undergone significant advances in the last few decades, see for example the textbooks  \cite{Toro:2001a}, \cite{Guinot:2012a} and \cite{Castro:2019a}.  A broad array of numerical frameworks have been proposed over the years, including finite difference methods, finite volume methods and finite element methods. More refined distinctions include centred methods and upwind methods.  One of the earliest, modern centred methods for solving the sallow water equations is due to Garc{\'i}a-Navarro et al. \cite{Garcia:1992a}, while examples of early modern upwind methods include References \cite{Glaister:1987a} and \cite{Toro:1992b}.  More recent works include References \cite{Vanzo2016} and \cite{Siviglia_2022}, the latter being close to the themes of the present paper.  A review of numerical methods for the shallow water equations is found in \cite{Toro:2007a}.\\

In the algorithm development process new research issues have arisen over the years. Source terms, due for example to chemical reactions and the bathymetry, have entered the research scene prominently, by posing important challenges to the numerical analyst. The concept of well-balanced scheme, in the presence of geometric-type source terms, has become a prominent topic. As a matter of fact, it is a mandatory requirement for a numerical method intended for solving the shallow water equations these days.  Key contributions in this topic include \cite{Bermudez:1994a}, \cite{Vazquez:1999a}, \cite{Pares:2004a},  \cite{Berthon:2016a},
\cite{castro2017well},
\cite{castro2020well}
 and many more. The concept of upwinding, classical for determining the numerical flux, was first introduced by Roe \cite{Roe:1986c} also for the source terms, see also \cite{Bermudez:1994a}.  Dealing numerically with dry fronts is another challenging topic for the design of numerical algorithms, see for example \cite{Castro:2006a}. The development of very high order numerical methods for the shallow water equations is yet another aspect of major importance, see for example \cite{Castro:2012a}. The introduction of adaptive mesh refinement (AMR) techniques in aerospace sciences \cite{Clarke:1993a} has also seen its importation into the shallow water territory,  see for example \cite{Beisiegel:2020}.  
Further algorithmic developments have been stimulated by the need for computationally efficient and flexible numerical strategies. Such development avenues are varied and include, for example, the use of local timestepping \cite{Dazzi2018}, \cite{Dumbser:2007c} as well as modification/reformulation of the governing equations via acceleration factors \cite{Carraro2018} or hyperbolization techniques \cite{Toro:2014a}, \cite{Vanzo2016}, to mention but a few.  A further strategy is furnished by explicit-implicit approaches, whereby the disparity in wave speeds can be exploited to allow for larger time steps, leading to gains in efficiency.  Potentially,  the present work could be exploited in this framework, as anticipated in \cite{Dumbser:2016d} for the compressible Navier-Stokes equations.  See also \cite{Boscheri:2021a} and \cite{Boscheri:2021b}.

In this paper we adopt the upwind philosophy for designing numerical methods for the shallow water equations. More specifically, we adopt the flux vector splitting approach (FVS) first proposed by Toro and Vazquez \cite{Toro:2012b} for the Euler equations of gas dynamics. It is important to realize that this is not the time-splitting method, which is a popular methodology for dealing with sources or multiple space dimensions, an example of which is the Strang splitting. See Chap. 15 of  \cite{Toro:2009a} for details on the time splitting method.  For background on FVS methods see Chap. 8 of \cite{Toro:2009a}. The particular flux splitting method proposed in \cite{Toro:2012b}, called hereafter the TV splitting,  has the distinctive feature of separating all pressure terms from the advection terms in the equations. When deploying such approach to the shallow water equations, however, two modifications are introduced in this paper. The first regards the continuity equation, which in the original TV splitting was kept in the advection system. In this paper such term is moved to the pressure system. This is justified in the framework of deriving zero-dimensional models from one-dimensional models, a popular approach in computational haemodynamics \cite{Formaggia:2009a}, \cite{Mueller:2014a}, \cite{Toro:2021a}.  The zero-dimensional models derived in this manner keep the advection term of the continuity equation and the advection operator vanishes when the inertial term in the momentum equation in the full system is discarded. Moreover, the straight extension of the TV splitting, so good for the Euler equations, does not work well for the shallow water equations, anomaly which is actually remedied by the proposed modification. The second aspect that needs a modification of the original TV splitting concerns non-conservative forms of the equations; these appear naturally when extending the shallow water equation to include, for example, sediment transport \cite{Siviglia_2022}. In this case there is no flux to split and therefore the scheme becomes a splitting of advection terms from pressure terms.  This second aspect is not relevant to the contents of the present paper, as the equations adopted here have conservation-law form.  Given these modifications we shall, generally,  speak of  TV-type splitting schemes.  \\

The present  paper is structured as follows: In Sect. 2, the shallow water equations are first briefly reviewed; then the basic new splitting scheme is formulated for the one-dimensional shallow water equations augmented by an advection equation for a generic passive scalar. The $3 \times 3$ system of equations is split into two systems of partial differential equations, namely an advection system and a pressure system.  For conservative systems, the sought,  single intercell numerical flux  is found through contributions from the advection system and the pressure system.  In Sect. 3 we provide solutions for the simplified pressure system at the element interface, which in turns provides the advection velocity for the advection system at such interface, completing the sought numerical flux. The resulting scheme is probably the simplest, complete upwinding scheme available in the literature; that is,  upwinding is provided by the scheme for every characteristic field in the equations, without the complexity of the Godunov method with the exact Riemann solver \cite{Toro:2001a}, or some other approximate,  non-linear complete Riemann solver. In Sect. 4 we extend the scheme to the two-dimensional shallow water equations with bathymetry source terms, on unstructured triangular meshes; we also extend the full methodology to high-order of accuracy in space and time following the ADER approach \cite{Toro:2001c}. Schemes of up to fifth order of accuracy in space and time are constructed; a convergence rates study verifies that such orders of accuracy are actually attained satisfactorily. In Sect. 5 we assess the methods on a suite of test problems with reference solution, including Riemann problems and a test involving trans-critical flow over a bed bump, with a stationary shock over the bump. Numerical results from the splitting scheme of this paper are compared against exact solutions and against the Godunov upwind method used in conjunction with the exact Riemann solver. In Sect. 6 we first present numerical results from the high-order schemes on two-dimensional test problems with reference solution. The section is finalised with an application of the full methodology to a tsunami wave propagation in the Pacific Ocean, in a realistic scenario.  Conclusions are drawn in Sect. 7.

\section{The 1D shallow water equations}

In this section we present the flux-vector splitting scheme as applied to the classical  shallow water equations, in one spatial dimension,  augmented by an advection equation for a passive scalar.  Sect. \ref{sec:PressureSystem} completes the presentation for the 1D case, while the two-dimensional extension is presented in Sect.  \ref{sec:2D}.

\subsection{Governing equations and the novel flux splitting}

\begin{figure}[htb]
\centering
\includegraphics[width=0.5\columnwidth]{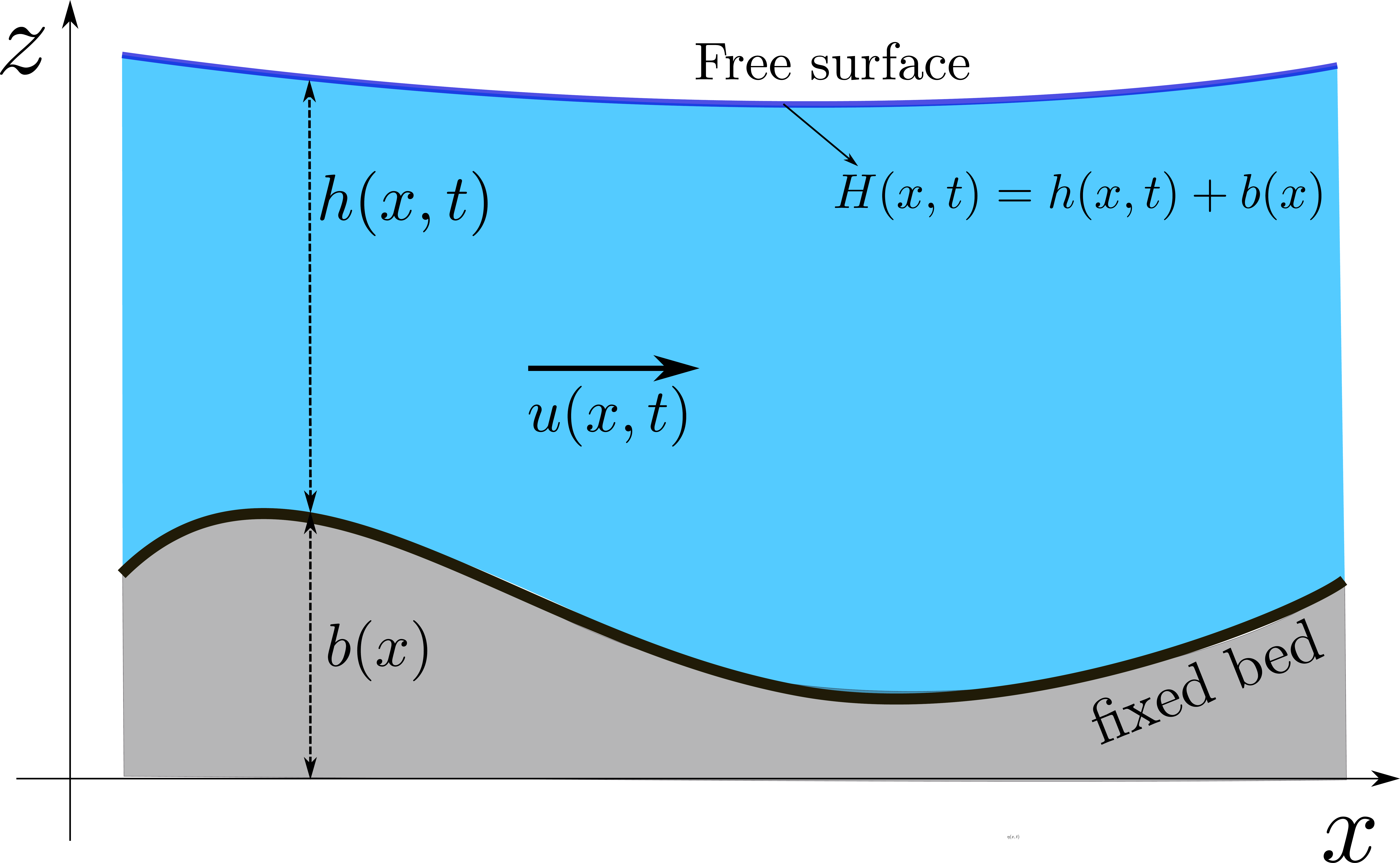}
\caption{Schematic of the  free-surface water flow problem over a prescribed fixed, variable bottom represented by $b(x)$. 
}
\label{fig:Sketch}
\end{figure}

The 1D governing equations are obtained under the well-know shallow water conditions \cite{Toro2001book}, 
\cite{Guinot:2012a}, \cite{Castro:2019a}.  Consider a Cartesian reference system $(x,z)$,  in which the $z$ [\si{m}] axis is vertical, while the horizontal axis $x$ [\si{m}] represents the longitudinal coordinate; see Fig.~\ref{fig:Sketch}. The governing equations include  conservation of water mass (continuity equation)
\begin{equation}                         \label{eq:1D1}
      \partial_{t}h +\partial_{x}(hu)=   0 \;
\end{equation}
and momentum balance
\begin{equation}                         \label{eq:1D2}
      \partial_{t}(hu) + \partial_{x}( hu^{2}+ \frac{1}{2}gh^{2})= -gh \partial_{x}b \;, 
\end{equation} 
 where $t$ [\si{s}] is time and $g$=9.81 [\si{ms^{-2}}] is the acceleration due to gravity. The quantities involved are illustrated in Fig. \ref{fig:Sketch}. Here 
 $h$ [\si{m}] is the flow depth,  $b$ [\si{m}] is the bed elevation, $H(x,t) = h(x,t) + b(x)$ [\si{m}] is the free surface position and $u$ [\si{ms^{-1}}] is depth-averaged flow velocity.  Note that in the 2D case $b=b(x,y)$.  The flow discharge per unit width is defined as $q=uh$=[\si{m^2s^{-1}}]. In this work  friction terms are neglected.
We augment the classical equations (\ref{eq:1D1})-(\ref{eq:1D2}) with an additional equation for a passive scalar $\psi(x,t)$, advected with the fluid speed $u$. Written in conservation-law form this equation reads
\begin{equation}                         \label{eq:1D3}
      \partial_{t}(h\psi) + \partial_{x}(hu\psi)= 0 \;. 
\end{equation}
 For details see \cite{Toro:2001a}. The resulting $3 \times 3$ system may be written in vectorial form as
\begin{equation}                         \label{eq:1D4}
      \partial_{t}{\bf Q} +\partial_{x}{\bf F}({\bf Q})={\bf S}  \;, 
\end{equation}
with
\begin{equation}                         \label{eq:1D5}
       {\bf Q}=\begin{bmatrix}
       h \\ 
       hu\\
       h \psi
       \end{bmatrix} \;, \hspace{5mm} 
       {\bf F}({\bf Q})= \begin{bmatrix}
		hu \\
	    hu^{2}+\frac{1}{2}gh^{2}\\
		hu\psi
	\end{bmatrix} \;, \hspace{5mm} 
	{\bf S}({\bf Q})=\begin{bmatrix}
       0 \\ 
       -gh b'(x) \\
       0
       \end{bmatrix}
	   \;.
\end{equation}
In this part of the paper, we are primarily interested in the principal part of Eqs. (\ref{eq:1D4})-(\ref{eq:1D5}) and therefore we restrict ourselves to
the homogeneous case $\mathbf{S}$($\mathbf{Q}$) = ${\bf 0}$.

Following Toro and V\'azquez \cite{Toro:2012b} we apply their TV splitting, suitably modified,  to the flux vector ${\bf F}({\bf Q})$, as follows
\begin{equation}                         \label{eq:1D6}
       {\bf F}({\bf Q})= \begin{bmatrix}
		hu \\
	    hu^{2}+\frac{1}{2}gh^{2}\\
		hu\psi
	\end{bmatrix} 
       ={\bf F}^{(a)}+{\bf F}^{(p)}
	   \;,
\end{equation}
with 
\begin{equation}                         \label{eq:1D7}
	{\bf F}^{(a)}
	=\begin{bmatrix}
       0 \\ 
       hu^{2} \\
       hu\psi
       \end{bmatrix} \;, \hspace{3mm}
       {\bf F}^{(p)}
       =\begin{bmatrix}
       hu \\ 
       \frac{1}{2}gh^{2} \\
       0
       \end{bmatrix}
	   \;.
\end{equation}
Again, following \cite{Toro:2012b}, we write two coupled systems of equations as follows
\begin{equation}                         \label{eq:1D8}
      \left.
      \begin{array}{ccc}
      \mbox{i) advection system} &: & \partial_{t}{\bf Q} +\partial_{x}{\bf F}^{(a)}({\bf Q})={\bf 0}  \;, \\
      \\
        \mbox{ii) pressure system} &: &\partial_{t}{\bf Q} +\partial_{x}{\bf F}^{(p)}({\bf Q})={\bf 0} \;.
      \end{array} \right\}
\end{equation}
We then propose  solution strategies to determined numerical fluxes for each of the two split sub-systems \eqref{eq:1D8}, in order to obtain the numerical flux for the full unsplit system (\ref{eq:1D4}).

\subsection{The numerical scheme}
The numerical solution of the full system (\ref{eq:1D4}) is obtained through the conservative scheme
\begin{equation}                         \label{eq:1D9}
    {\bf Q}_{i}^{n+1} = {\bf Q}_{i}^{n}-\frac{\Delta t}{\Delta x}[ {\bf F}_{i+\frac{1}{2}}-{\bf F}_{i-\frac{1}{2}} ]  \;, 
\end{equation}
with numerical flux made up of the split fluxes, namely
\begin{equation}                         \label{eq:1D10}
    {\bf F}_{i+\frac{1}{2}}= {\bf F}_{i+\frac{1}{2}}^{(a)} + {\bf F}_{i+\frac{1}{2}}^{(p)}  \;. 
\end{equation}
The advection flux ${\bf F}_{i+\frac{1}{2}}^{(a)}$ and the pressure flux ${\bf F}_{i+\frac{1}{2}}^{(p)}$ are obtained from the advection and the pressure systems in (\ref{eq:1D8}), respectively. \\

For given intial conditions ${\bf Q}_{L}$, ${\bf Q}_{R}$ either side of the interface $x_{i+\frac{1}{2}}$ that naturally define a Riemann problem, see (\ref{eq:1D20}),  we propose the following formulation for the numerical fluxes in \eqref{eq:1D10}:
\begin{equation}                         \label{eq:1D11}
     {\bf F}_{i+\frac{1}{2}}^{(a)} = \left\{\begin{array}{c}
     \begin{bmatrix}
     0\\
      \displaystyle{q_{*} u_{L}  } \\
     q_{*} \psi_{L}
     \end{bmatrix}   \mbox{ if } q_{*} \ge 0 \;, \\
     \\
     \begin{bmatrix}
     0\\
     \displaystyle{q_{*}u_{R}  } \\
     q_{*}\psi_{R}
     \end{bmatrix}  \mbox{ if } q_{*} < 0 \;,
     \end{array} \right. 
\end{equation}
and
\begin{equation}                         \label{eq:1D12}
     {\bf F}_{i+\frac{1}{2}}^{(p)} =
     \begin{bmatrix}
     q_{*} \\
    \frac{1}{2}gh_{*}^{2}\\
    0
     \end{bmatrix}   \;.
\end{equation}
Here $h_{*}$ and $q_{*}$ are values along the $t$-axis obtained from solving the Riemann problem for the pressure system.  Such values are yet to be found. With such formulation, it turns out that the determination of the two split fluxes relies entirely on the solution of the Riemann problem for the pressure system alone.

\section{The pressure system and the Riemann problem}
\label{sec:PressureSystem}

\begin{figure}
      \centerline{
      \includegraphics[width=0.4\columnwidth]{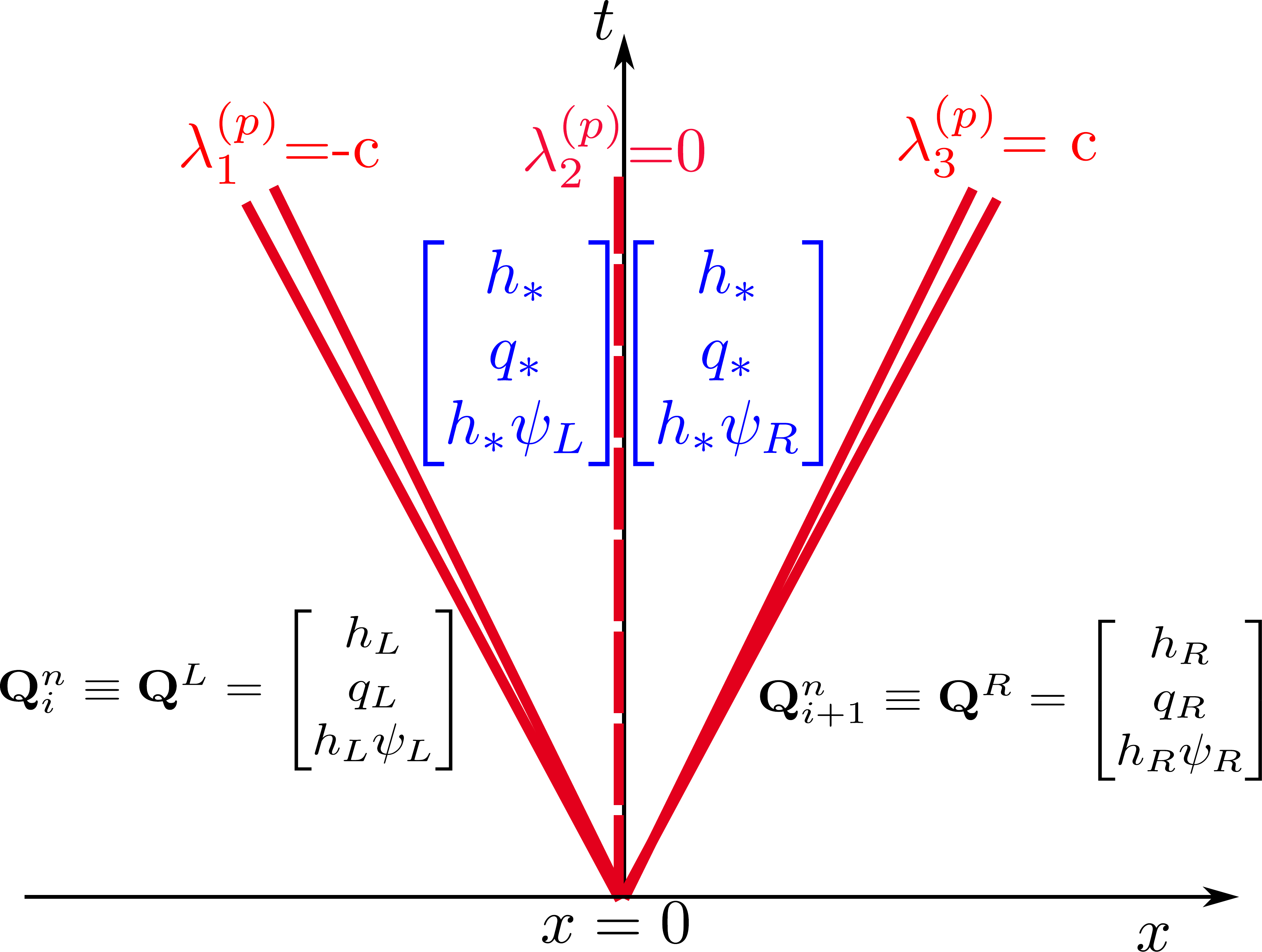}
       }
\caption{ \textbf{Structure of the solution of the Riemann problem for the pressure system resulting from the flux splitting}. There are two non-linear wave
families and a stationary contact discontinuity coinciding
with the $t$-axis. The wave pattern is always subcritical therefore the Godunov state for flux evaluation is always the star state, hence solution sampling is not required. 
}
\label{fig:RP}
\end{figure}

Consider the Riemann problem for the pressure system
\begin{equation} \label{eq:1D20}
\left.
\begin{array}{l}
       \partial_{t}{\bf Q} + \partial_{x} {\bf F}^{(p)}({\bf Q})={\bf 0} \;, \\
        \\
        {\bf Q}(x,0)=
     \left\{\begin{array}{ccc}
        {\bf Q}_{L}=\begin{bmatrix}
       h_{L} \\ 
       h_{L}u_{L}\\
       h_{L} \psi_{L}
       \end{bmatrix}  & \mbox{  if  } & x<0 \;, \\
       \\
        {\bf Q}_{R}=\begin{bmatrix}
       h_{R} \\ 
       h_{R}u_{R}\\
       h_{R} \psi_{R}
       \end{bmatrix} & \mbox{  if  }  & x>0 \;. \\
    \end{array} \right.
\end{array}\right\}
\end{equation}
The structure of the solution of  (\ref{eq:1D20}) is depicted in Fig.~\ref{fig:RP}. The wave pattern is always {\bf subcritical} and the sought {\bf Godunov state} ${\bf Q}_{i+\frac{1}{2}}(0)$ along the $t$-axis for flux evaluation is always identical to the intermediate {\bf Star State} between the two pressure waves, denoted 
here as ${\bf Q}_{*}=[h_{*}, q_{*}, \psi_{*} h_{*}]^{T}$.  Note here that in the {\bf Star Region} there are two values for $\psi_{*}$, which change discontinuously across the $t$-axis.  \\


In order to solve the Riemann problem (\ref{eq:1D20}) we need some preliminaries, such as the eigenstructure of the system. When written in quasi-linear form the pressure system becomes
\begin{equation}                         \label{eq:1D15}
       \partial_{t}{\bf Q} +{\bf J}^{(p)} \partial_{x}{\bf Q}={\bf 0} \;,  \hspace{3mm}
        {\bf J}^{(p)} = \begin{bmatrix}
		0   & 1 & 0\\
		gh & 0 & 0 \\
		0   & 0 & 0 
	\end{bmatrix} \;.
\end{equation}
The eigenvalues of the Jacobian matrix ${\bf J}^{(p)}$ are all real and given as
\begin{equation}                         \label{eq:1D16}
       \lambda_{1}^{(p)}=-c  \;, \hspace{3mm} \lambda_{2}^{(p)}=0  \;, \hspace{3mm} \lambda_{3}^{(p)}=c \;,
\end{equation}
where 
\begin{equation}                         \label{eq:1D17}
       c = \sqrt{g h}
\end{equation}
is the usual celerity in the shallow water equations. Note that the pressure system is always subcritical (subsonic), as seen in Fig.~\ref{fig:RP}, that is
\begin{equation}                         \label{eq:1D18}
       \lambda_{1}^{(p)}=-c < \lambda_{2}^{(p)}=0 <\lambda_{3}^{(p)}=c \;.
\end{equation}
The corresponding right eigenvectors are linearly independent and, after appropriate scaling,  are given as
\begin{equation}                         \label{eq:1D19}
       {\bf R}_{1}= \begin{bmatrix}
		1 \\
	    -c\\
		0
	\end{bmatrix} \;, \hspace{3mm}
	{\bf R}_{2}= \begin{bmatrix}
		0 \\
	    0\\
		1
		\end{bmatrix} \;, \hspace{3mm}
	{\bf R}_{3}= \begin{bmatrix}
		1 \\
	    c \\
		0
		\end{bmatrix} \;.
\end{equation}
The pressure system (\ref{eq:1D15}) is therefore strictly hyperbolic.

In what follows  we propose two approaches to solve the Riemann problem for the pressure system,  to find  suitable expressions for $h_{*}$ and  $q_*$. These are the primary variables, as the solution for $\psi_{*}$ is trivial.

\subsection{Simple waves and generalised Riemann invariants}

Assuming the solution of the Riemann problem consists of smooth waves, we apply generalised Riemann invariants.  For $\lambda_{1}^p=-c$ we have
\begin{equation}                         \label{eq:1D21}
       \frac{2}{3}h \sqrt{gh} + q =\mbox{constant} \;,  \hspace{3mm} \psi = \mbox{constant}.
\end{equation}
This allows us to connect the left state ${\bf Q}_{L}$ to the unknown state ${\bf Q}_{*}$, namely
\begin{equation}                         \label{eq:1D22}
       \frac{2}{3}h_{*} \sqrt{gh_{*}} + q_{*} = \frac{2}{3}h_{L} \sqrt{gh_{L}} + q_{L}  \;,
       \hspace{3mm} \psi_{*L} =\psi_{L} \;.
\end{equation}
From here we may write
\begin{equation}                         \label{eq:1D23}
       q_{*} =  q_{L}  - f_{L}(h_{*}, h_{L}) \;, \hspace{4mm} f_{L}(h_{*}, h_{L})=\frac{2}{3}\sqrt{g}(h_{*}^{\frac{3}{2}}-h_{L}^{\frac{3}{2}}) \;,
\end{equation}
an expression that will be useful when solving the Riemann problem exactly (iteratively).

The generalised Riemann invariant corresponding to $\lambda_{3}^p=c$ gives 
\begin{equation}                         \label{eq:1D24}
       \frac{2}{3}h \sqrt{gh} - q = \mbox{constant} \;,  \hspace{3mm} \psi = \mbox{constant}.
\end{equation}
This is used to connect the right state ${\bf Q}_{R}$ to the unknown state ${\bf Q}_{*}$; it gives
\begin{equation}                         \label{eq:1D25}
       \frac{2}{3}h_{*} \sqrt{gh_{*}} - q_{*} = \frac{2}{3}h_{R} \sqrt{gh_{R}} - q_{R}  \;, 
        \hspace{3mm} \psi_{*R} =\psi_{R} \;.
\end{equation}
From here we write
\begin{equation}                         \label{eq:1D26}
       q_{*} =  q_{R}  + f_{R}(h_{*}, h_{R}) \;, \hspace{4mm} f_{R}(h_{*}, h_{R})=\frac{2}{3}\sqrt{g}(h_{*}^{\frac{3}{2}}-h_{R}^{\frac{3}{2}})\;.
\end{equation}

The exact solution of the Riemann problem will also include shock waves and, in general,  will be iterative. However, assuming the special case in which both pressure waves are rarefaction waves, we can obtain a closed-form, approximate solution.  Collecting (\ref{eq:1D22}) and (\ref{eq:1D25})  we have 
an algebraic system for two unknowns, namely
\begin{equation}                         \label{eq:1D28}
      \left.\begin{array}{l}
       q_{*} + \frac{2}{3}h_{*} \sqrt{gh_{*}}= q_{L} + \frac{2}{3}h_{L} \sqrt{gh_{L}}\;,\\
       \\
       q_{*} - \frac{2}{3}h_{*} \sqrt{gh_{*}}= q_{R} - \frac{2}{3}h_{R} \sqrt{gh_{R}}
       \;.
       \end{array}\right\}
\end{equation}
The exact solution of (\ref{eq:1D28})  is
\begin{equation}                         \label{eq:1D29}
 \left.
   \begin{array}{l}
   q_{*} = \frac{1}{2}(q_{L}+q_{R}) + \frac{1}{3} \sqrt{g}(h_{L}^{\frac{3}{2}}-h_{R}^{\frac{3}{2}})  \;,\\
   \\
   h_{*} = [\frac{1}{2}(h_{L}^{\frac{3}{2}}+h_{R}^{\frac{3}{2}}) -\frac{3}{4\sqrt{g}}(q_{R}-q_{L}) ]^{\frac{2}{3}}\;.
   \end{array}
\right\}
\end{equation}
This solution is called the {\bf two-rarefaction solution}. As a matter of fact this approximate solution is quite accurate also in the case in which the solution of the Riemann problem involves shock waves.
It is recommended for practical use.

\subsection{Shocks and Rankine-Hugoniot conditions}

For the exact solution of the Riemann problem we also need to consider the possible presence of shock waves in the wave pattern that emerges from the solution. The final solution however is not direct, it is iterative, as we shall see.\\

\paragraph{Left-facing shock wave}

For a left shock wave of speed $S_{L}$ connecting ${\bf Q}_{L}$ to the solution in the star region ${\bf Q}_{*}$, afer applying Rankine-Hugoniot conditions, we may write
\begin{equation}                         \label{eq:1D30}
       q_{*} - q_{L} =S_{L} (h_{*} - h_{L}) \;, \hspace{4mm}  \frac{1}{2}g(h_{*}^{2} - h_{L}^{2})=S_{L} (q_{*} - q_{L}) \;.
\end{equation}
From here we obtain
\begin{equation}                         \label{eq:1D31}
       S_{L}^{2}= \frac{1}{2}g(h_{*} + h_{L})  \;, 
\end{equation}
from which 
\begin{equation}                         \label{eq:1D31a}
       S_{L}= \pm \sqrt{\frac{1}{2}g(h_{*} + h_{L})} = \pm \sqrt{\frac{1}{2}(c_{*}^{2} + c_{L}^{2})} 
       =  \pm c_{L} \sqrt{\frac{1}{2}\left(\frac{h_{*}}{h_{L}} +1\right) }\;.
\end{equation}
The choice of the correct sign in (\ref{eq:1D31a}) is determined from invoking the entropy condition, which for a left entropy-satisfying shock states
\begin{equation}                         \label{eq:1D31b}
       \lambda_{1}({\bf Q}_{L}) > S_{L}> \lambda_{1}({\bf Q}_{*}) \;.
\end{equation}
That is $-c_{L} > S_{L}>- c_{*}$,  from which it follows
\begin{equation}                         \label{eq:1D31c}
           S_{L}<0            \mbox{     and      }         h_{*}> h_{L} \;. 
 \end{equation}
Therefore  the correct choice of sign in (\ref{eq:1D31a}) is minus. That is
\begin{equation}                         \label{eq:1D32}
       S_{L}= -\sqrt{\frac{1}{2}g(h_{*} + h_{L})} \;.
\end{equation}
In readiness for solving the Riemann problem, from (\ref{eq:1D30}) and (\ref{eq:1D32}) we write the expression
\begin{equation}                         \label{eq:1D33}
       q_{*} = q_{L} - f_{L} (h_{*},h_{L}) \;, \hspace{4mm}  
       f_{L} (h_{*},h_{L}) = \sqrt{\frac{1}{2}g(h_{*}+ h_{L})}(h_{*} - h_{L}) \;
\end{equation}
to be used later.

\paragraph{Right-facing shock wave}
Analogously, for a right shock wave of speed $S_{R}$ connecting ${\bf Q}_{R}$ to the solution in the star region with state ${\bf Q}_{*}$ we have
\begin{equation}                         \label{eq:1D34}
       S_{R}^{2}= \frac{1}{2}g(h_{*} + h_{R}) \;,
\end{equation}
leading to
\begin{equation}                         \label{eq:1D34a}
       S_{R}= \sqrt{\frac{1}{2}g(h_{*} + h_{R})} =  \sqrt{\frac{1}{2}(c_{*}^{2} + c_{R}^{2})} 
       =  c_{R} \sqrt{\frac{1}{2}\left(\frac{h_{*}}{h_{R}} +1\right) }\;.
\end{equation}
Again, in readiness for solving the Riemann problem exactly we write the expression
\begin{equation}                         \label{eq:1D36}
       q_{*} = q_{R} + f_{R} (h_{*},h_{R}) \;, \hspace{4mm}  
       f_{R} (h_{*},h_{R}) = \displaystyle{ \sqrt{\frac{1}{2}g(h_{*}+ h_{R})}(h_{*} - h_{R}) } \;.
\end{equation}

\subsection{Exact iterative solution of the Riemann problem}

From equations (\ref{eq:1D23}), (\ref{eq:1D26}), (\ref{eq:1D33}) and (\ref{eq:1D36}) we may write a single non-linear
algebraic equation for $h=h_{*}$, namely
\begin{equation}                         \label{eq:1D37}
       f(h) = f_{L} (h,h_{L}) +f_{R} (h,h_{R}) +q_{R}-q_{L}=0 \;, 
\end{equation}
where
\begin{equation}                         \label{eq:1D38}
       f_{L} (h,h_{L}) =\left\{\begin{array}{ccc}
         \frac{2}{3}\sqrt{g}(h_{*}^{\frac{3}{2}}-h_{L}^{\frac{3}{2}}) &  \mbox{  if  } & h \le h_{L}\\
          \\
          \sqrt{\frac{1}{2}g(h_{*}+ h_{L})}(h_{*} - h_{L})               &  \mbox{  if  } & h > h_{L}
          \end{array}\right.
\end{equation}
and
\begin{equation}                         \label{eq:1D39}
       f_{R} (h,h_{R}) =\left\{\begin{array}{ccc}
         \frac{2}{3}\sqrt{g}(h_{*}^{\frac{3}{2}}-h_{R}^{\frac{3}{2}}) &  \mbox{  if  } & h \le h_{R}\\
          \\
          \sqrt{\frac{1}{2}g(h_{*}+ h_{R})}(h_{*} - h_{R})                &  \mbox{  if  } & h > h_{R}
          \end{array}\right.
\end{equation}
Equation (\ref{eq:1D37}) may be solved by a Newton-Raphson method, that is
\begin{equation}                         \label{eq:1D40}
       h^{(k+1)} = h^{(k)}  - \frac{f(h^{(k)}  )}{f'(h^{(k)}  )}\;, 
\end{equation}
where the derivatives are given as
\begin{equation}                         \label{eq:1D41}
       f_{K}'= \left\{  \begin{array}{ccc} 
       \sqrt{gh}   & \mbox{  if  } & h \le h_{K} \;, \\
       \\
       \displaystyle{\sqrt{\frac{g}{8}}  \frac{(3h+ h_{K})}{\sqrt{h+ h_{K}}}} & \mbox{  if  } & h > h_{K} \;.
       \end{array}\right.
\end{equation}
Here $h_{K}$ is either $h_{L}$  or $h_{R}$ .
As a guess value $h^{(0)}$ to start iteration (\ref{eq:1D40}) we use the two-rarefaction solution  (\ref{eq:1D29}). The iteration is stopped when the following condition is met
\begin{equation}                         \label{eq:1D42}
          \Delta h = \frac{|h^{(k+1)} - h^{(k)}|}{\frac{1}{2}(h^{(k)} + h^{(k+1)})} \le TOL \;,
\end{equation}
where $TOL$ is a prescribed tolerance, e.g.  $TOL=10^{-9}$.\\

Once the solution for $h_{*}$ has been obtained from solving   (\ref{eq:1D37}) iteratively, the solution for 
$q_{*}$ is obtained from  (\ref{eq:1D23}),  (\ref{eq:1D26}), (\ref{eq:1D33}) and (\ref{eq:1D36})  as
\begin{equation}                         \label{eq:1Dqstar}
       q_{*} = \frac{1}{2}(q_{L}+  q_{R})+  \frac{1}{2}[f_{R} (h_{*},h_{R}) - f_{L} (h_{*},h_{L})] \;.
\end{equation}

\noindent{\bf Remarks:}

\begin{enumerate}

\item We have provided two ways to calculate the solution of the Riemann problem and provide the sought Godunov state  for use in the pressure and advection fluxes to make up the full numerical flux.  One option is the direct solution (\ref{eq:1D29}) and another is the iterative solution from solving  (\ref{eq:1D37}), followed by the direct calculation of $q_{*}$ from (\ref{eq:1Dqstar}).

\item The practical implementation of both methodologies reveal that  there is no visible difference between their respective  numerical results.  Therefore we suggest to use the approximate solution.

\item Note the value for $\psi_{*}$, which is either $\psi_{L}$ or $\psi_{R}$. See Fig.  \ref{fig:RP}.

\item For the two-dimensional case, the transverse component of velocity, say $v$, behaves identically as $\psi$ in the solution of the Riemann problem.

\end{enumerate}

\section{Extension to two dimensions and high order of accuracy}
\label{sec:2D}

In this section  the proposed TV-type splitting scheme for the shallow water equations  is extended to solve the two-dimensional equations on unstructured meshes, to first order and to high-order order of accuracy in space and time.  

\subsection{The equations in two space dimensions}

The equations are formulated in a Cartesian reference system $(x,y,z)$, where the $z$ axis is vertical and the $(x,y)$ plane is horizontal. The equation for conservation of fluid mass reads 
\begin{equation}\label{eq:2D_gov_eq_mass}
\partial_{t}h +\partial_{x}(q_x)+\partial_{y}(q_y)=   0, \;
\end{equation}
whilst the equations for balance of momenta in $x$ and $y$ directions read:  
\begin{equation}\label{eq:2D_gov_eq_momx}
\partial_{t}(q_x) + \partial_{x}\left( \dfrac{{q_x}^2}{h} + \frac{1}{2}gh^{2} \right) +\partial_{y}\left(\dfrac{{q_x}{q_y}}{h}\right) = -gh\partial_{x} b \;, 
\end{equation}
and
\begin{equation}\label{eq:2D_gov_eq_momy}
\partial_{t}(q_y) +\partial_{x}\left(\dfrac{{q_x}{q_y}}{h}\right) + \partial_{y}\left( \dfrac{{q_y}^2}{h} + \frac{1}{2}gh^{2} \right)  = -gh\partial_{y} b \;.
\end{equation}
Here $h(x,y,t)$ and $b(x,y)$ [m] denote the water depth and bottom elevation respectively, whilst $q_x$ and $q_y$ are the directional components of the vector of specific discharge ${\mathbf q}=(q_x, q_y)$ [\si{m^2s^{-1}}]. The depth-averaged flow velocity vector can thus be expressed as ${\mathbf u}:=(u, v)=(q_x/h,q_y/h)$ [\si{ms^{-1}}]. The unknowns of the problem are the water depth $h(x,y,t)$ and the discharge per unit width (or momentum), that is its components $q_x(x,y,t)$ and $q_y(x,y,t)$, from which the the velocity components  $u$ and   $v$ are calculated. The bed elevation $b(x,y)$ is prescribed, assumed arbitrary and fixed.  The resulting $3 \times 3$ system of governing equations \eqref{eq:2D_gov_eq_mass}, \eqref{eq:2D_gov_eq_momx} and \eqref{eq:2D_gov_eq_momy}, can be written in vectorial form as
\begin{equation}\label{eq:2D_vectorial}
\partial_{t}{\bf Q} + \partial_{x}{\bf F}_x({\bf Q}) + \partial_{y}{\bf F}_y({\bf Q}) ={\bf S}  \;, 
\end{equation}
with
\begin{equation}\label{eq:2D_vec_terms}
\renewcommand{\arraystretch}{1.5}
{\bf Q}=\begin{bmatrix}
h \\ 
q_x\\
q_y\\
\end{bmatrix} \;, \hspace{5mm} 
{\bf F}_{x}({\bf Q})= \begin{bmatrix}
q_x\\
\dfrac{{q_x}^2}{h} + \frac{1}{2}gh^{2}\\
\dfrac{{q_x}{q_y}}{h}\\
\end{bmatrix} \;, \hspace{5mm}
{\bf F}_{y}({\bf Q})= \begin{bmatrix}
q_y \\
\dfrac{{q_x}{q_y}}{h}\\
\dfrac{{q_y}^2}{h} + \frac{1}{2}gh^{2}\\
\end{bmatrix} \;, \hspace{5mm}
{\bf S}({\bf Q})=\begin{bmatrix}
0 \\ 
-gh\partial_{x}b \\
-gh\partial_{y}b \\
\end{bmatrix}
\;.
\end{equation}
The TV-type splitting procedure for the 1D case is applied directly by splitting the fluxes  into  \textit{advection} and \textit{pressure} contributions,  as follows
\begin{equation}\label{eq:2D_split1}
{\bf F}_x={\bf F}_x^{(a)}+{\bf F}_x^{(p)}\qquad,\qquad
{\bf F}_y={\bf F}_y^{(a)}+{\bf F}_y^{(p)}\;,
\end{equation}
with 
\begin{equation}\label{eq:2D_split2}
\renewcommand{\arraystretch}{2}
{\bf F}_x^{(a)}
=\begin{bmatrix}
0\\
\dfrac{{q_x}^2}{h}\\
\dfrac{{q_x}{q_y}}{h}\\
\end{bmatrix} \;, \qquad
{\bf F}_x^{(p)}
=\begin{bmatrix}
q_x\\
\frac{1}{2}gh^{2}\\
0\\
\end{bmatrix}
\;,\qquad
{\bf F}_y^{(a)}
=\begin{bmatrix}
0\\
\dfrac{{q_x}{q_y}}{h}\\
\dfrac{{q_y}^2}{h}\\
\end{bmatrix} \;, \qquad
{\bf F}_y^{(p)}
=\begin{bmatrix}
q_y\\
0 \\
\frac{1}{2}gh^{2}\\
\end{bmatrix}.
\end{equation}

\subsection{The 2D numerical scheme}

\begin{figure}[tbp]
    \centering
    \includegraphics[angle=0,width=0.5\columnwidth]{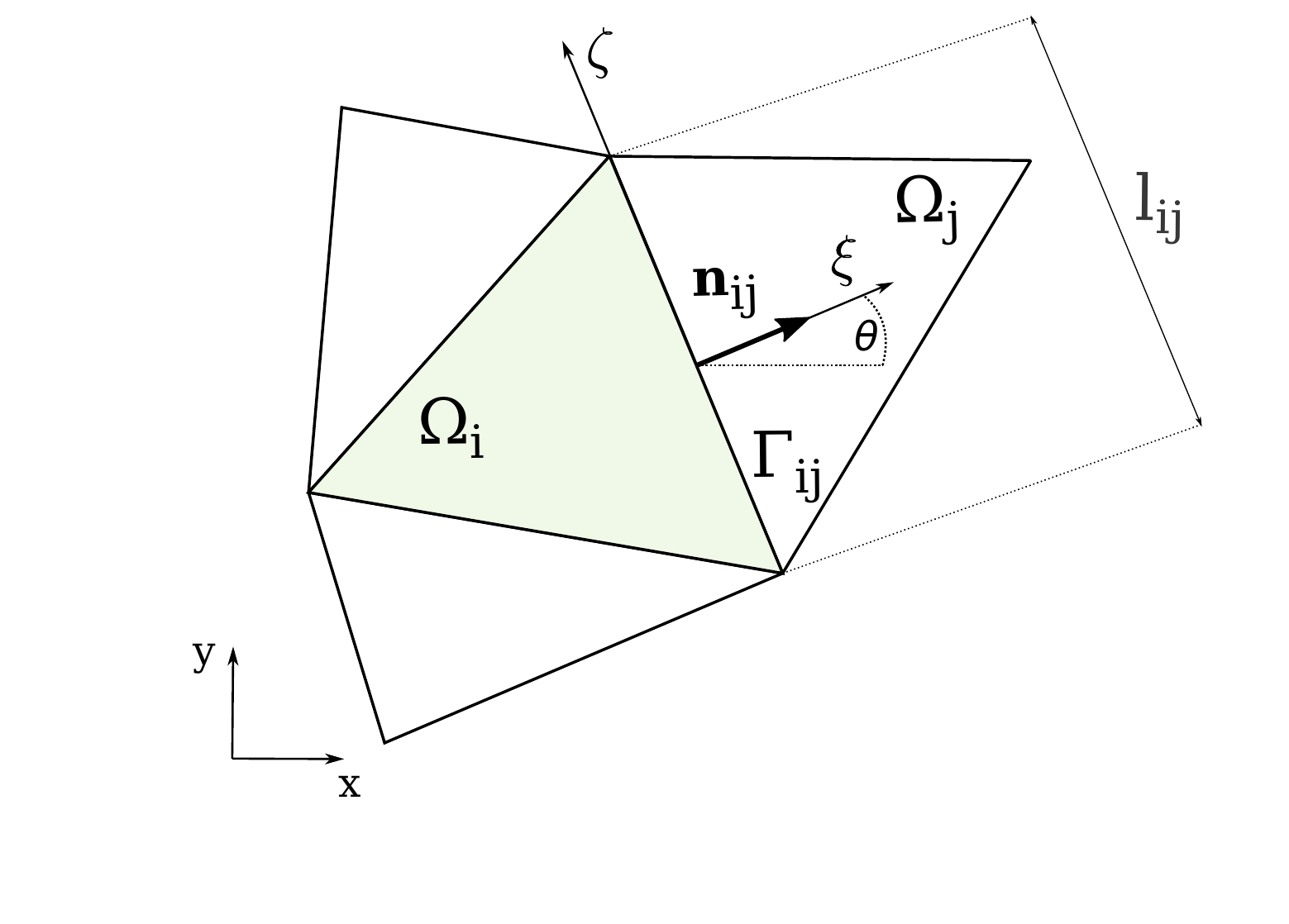}
    \caption{Sketch of the triangular discretization,  including notation for the computational cell $i$ and its neighbour $j$.}
    \label{fig:unstructured_discretization}
\end{figure}

The two dimensional domain of interest is discretised via an unstructured triangular mesh.  A conforming triangulation $T_{\Omega}$ of the computational domain $\Omega \subset \mathbb{R}^2$ by triangular elements $\Omega_i$ is assumed, such that $T_{\Omega}=\bigcup \Omega_i$; see Fig.~\ref{fig:unstructured_discretization}.  For a given finite volume triangular element $\Omega_i$,   $j=1, 2, 3$ denote the indexes of the neighbouring cells $\Omega_j$; $\Gamma_{ij}$ is the common edge between $\Omega_i$ and its neighbour $\Omega_j$ and $l_{ij}$ is its length. ${\mathbf{n}}_{ij}=(n_{ij,x},n_{ij,y})=(cos\theta, sin\theta)$ is the unit vector which is normal to the edge $\Gamma_{ij}$ and points toward the cell $\Omega_j$.  For each cell edge $\theta$ is the positive angle formed by the $x$-direction in the Cartesian mesh and the outward unit normal ${\mathbf{n}}_{ij}$. Choosing the $x$-direction as the reference direction is a convention. One could also choose the $y$-direction.

The general finite volume update formula applied to a general triangular element $\Omega_i$ reads:
\begin{equation}\label{eq:2D_update_formula}
{\bf Q}_i^{n+1}={\bf Q}_i^n-\frac{\Delta t}{\left|\Omega_i\right|}\sum_{j=1}^3 \int_{l_{ij}}\left[ \frac{1}{\Delta t} \int_{t^n}^{t^{n+1}} {\bf F}_{ij} \,\textnormal{d}t 
\right]\textnormal{d\textbf{x}} +\Delta t  \boldsymbol{S}_{i}\;.
\end{equation}
Here ${\bf F}_{ij}$ is the numerical flux normal to the cell interface $\Gamma_{ij}$; $\left|\Omega_i\right|$ is the area of $\Omega_i$ and $\boldsymbol{S}_{i}$ is the numerical source within the volume $\Omega_i$. In what follows we shall provide all the necessary ingredients to implement and apply scheme (\ref{eq:2D_update_formula}) in first-order mode. This is essentially the Godunov upwind method in two space dimensions, in which the exact solution of the Riemann problem is replaced by the TV-type split fluxes just described. The ADER higher-order extension of the scheme will also make use of the same one-step formula (\ref{eq:2D_update_formula}).

The numerical fluxes ${\bf F}_{ij}$ in \eqref{eq:2D_update_formula} are computed by exploiting the rotational invariance of the 2D  governing shallow water equations, see \cite{Toro:2001a}. In this manner the Riemann problem defined by initial piecewise constant data across the interface $\Gamma_{ij}$ can be posed in a  locally rotated reference system  $(\xi,\zeta)$ by means of a rotation ${\bf T}$ applied to the initial data; see Fig.~\ref{fig:unstructured_discretization}.  The final flux in the original reference frame ($x,y$) is obtained by means of an inverse rotation ${\bf T}^{-1}$ applied to the flux evaluated on the rotated data. That is
\begin{equation}\label{eq:rot_inv}
\bF_{ij}={\bf n}_{ij} \cdot [\bF_x(\bQ),\bF_y(\bQ)]={\bf T}^{-1}{\bF}_x({\bf T}(\bQ))\;.
\end{equation}
The rotation matrix ${\bf T}={\bf T}(\theta)$ and its inverse ${\bf T}(\theta)^{-1}$ are
\begin{equation}\label{eq:rotation_matrix}
{\bf T}=\begin{bmatrix}
1 & 0 & 0\\
0 & cos\theta & sin\theta\\
0 & -sin\theta & cos\theta\\
\end{bmatrix} \;, \hspace{3mm}
{\bf T}^{-1}=\begin{bmatrix}
1 & 0 & 0\\
0 & cos\theta & -sin\theta\\
0 & sin\theta & cos\theta\\
\end{bmatrix} \;.
\end{equation}
For details see  Sect. 3.10 of \cite{Toro:2001a}.\\

In practice we denote the locally rotated unknowns as $\bQ^r={\bf T}(\bQ)$ and the sought corresponding normal fluxes in the rotated system as $\bF_{ij}^{r}={\bF}_x({\bf T}(\bQ))$, namely
\begin{equation}\label{eq:rotated_var}
\renewcommand{\arraystretch}{2}
{\bQ^r}=\begin{bmatrix}
h \\ 
q_\xi\\
q_\zeta\\
\end{bmatrix}\;, \hspace{3mm}
\bF_{ij}^{r}= \begin{bmatrix}
q_\xi\\
\dfrac{{q_\xi}^2}{h} + \frac{1}{2}gh^{2}\\
\dfrac{{q_\xi}{q_\zeta}}{h}\\
\end{bmatrix}.
\end{equation}
 A solution of the generalized Riemann problem is then sought for the rotated interface fluxes $\bF_{ij}^{r}$. In analogy with 1D formulation and by definition \eqref{eq:rotated_var}, the rotated fluxes can be split into advection $(a)$ and pressure $(p)$ contributions as
\begin{equation}\label{eq:2D_num_split_flux}
\bF_{ij}^r= {\bF_{ij}^r}^{(a)} + {\bF_{ij}^r}^{(p)}.
\end{equation}
Given $h_{*}$ and $q_{\xi *}$, the solutions in the star region of the local Riemann problem associated with the pressure system, the final sought fluxes terms of \eqref{eq:2D_num_split_flux} can be expressed as
\begin{equation}\label{eq:2D_flux_adv_sol}
\renewcommand{\arraystretch}{2}
{\bF_{ij}^r}^{(a)} = \left\{\begin{array}{c}
\begin{bmatrix}
0\\
\dfrac{q_{\xi *}q_{\xi i}}{h_{i}} \\
\dfrac{q_{\xi *}q_{\zeta i}}{h_{i}} \\
\end{bmatrix}   \mbox{ if } q_{\xi *} \ge 0 \;, \\
\\
\begin{bmatrix}
0\\
\dfrac{q_{\xi *}q_{\xi j}}{h_{j}} \\
\dfrac{q_{\xi *}q_{\zeta j}}{h_{j}} \\
\end{bmatrix}   \mbox{ if } q_{\xi *} < 0,
\end{array} \right. 
\end{equation}
and
\begin{equation}\label{eq:2D_flux_press_sol}
{\bF_{ij}^r}^{(p)} =
\begin{bmatrix}
q_{\xi *} \\
\frac{1}{2}gh_*^{2}\\
0\\
\end{bmatrix}   \;.
\end{equation}
The procedure to evaluate the star-region values $h_{*}$ and $q_{\xi *}$ is analogous to the 1D formulation. Omitting the explicit derivation, the solution of the two-rarefaction solution for the 2D case (equivalent to \eqref{eq:1D29} for the 1D case), reads 
\begin{equation}\label{eq:2D_two_rar}
\left.
\begin{array}{l}
q_{*} = \frac{1}{2}(q_{\xi i}+q_{\xi j}) + \frac{1}{3} \sqrt{g}(h_i^{\frac{3}{2}}-h_j^{\frac{3}{2}})  \;,\\
\\
h_{*} = [\frac{1}{2}(h_{i}^{\frac{3}{2}}+h_{j}^{\frac{3}{2}}) -\frac{3}{4\sqrt{g}}(q_{\xi j}-q_{\xi i}) ]^{\frac{2}{3}}\;.
\end{array}
\right\}
\end{equation}
%


\subsection{High-order extension following the ADER approach}

So far we have described the method on triangular meshes in first-order mode. The  extension to high-order of accuracy is realised here via the ADER approach, in which the time updating one-step formula remains, identically,  that for the first-order scheme (\ref{eq:2D_update_formula}) The difference between the first-order scheme and any of the high-order extensions resides entirely in the numerical fluxes and the numerical source.  In the fully discrete ADER approach, it is possible to increase the space and time accuracy to any desired order. In other words, in the ADER approach there is no theoretical barrier to the order of accuracy and therefore one speaks of schemes of {\bf arbitrary accuracy}.  After some,  partial preliminary communications during the 1990s,  the ADER framework was comprehensively presented in \cite{Toro:2001c} for linear equations on Cartesian meshes in one and multiple space dimensions.  Results of up to $10$th order of accuracy in space and time were presented.  A historical note, all the material in \cite{Toro:2001c} was presented in 1999 in a plenary presentation at an international conference to celebrate Godunov's $70$th birthday, at the University of Oxford, UK.  Many contributions to the development and application of ADER have been reported since then.  To mention just a few examples, see
\cite{Toro:2002a},
\cite{Titarev:2002a},
\cite{Dumbser:2005a},
\cite{Dumbser:2008a},
\cite{Dumbser:2008b},
\cite{Dumbser:2016c},
 \cite{Toro:2015b},
 \cite{Dematte:2020a},
 \cite{Goetz:2016a},
 \cite{Goetz:2016b}.  An introduction to ADER is provided  in Chaps. 19 and 20 of \cite{Toro:2009a}.  For more recent  advances on  ADER schemes see the review \cite{Toro:2020a} and references therein.  \\
 
The high order ADER extension relies on the first-order scheme just described and on two new elements, namely (i) high-order non-linear spatial reconstruction and (ii)  the solution of the generalized Riemann problem at the cell interface.  For the spatial reconstruction we make use of a  WENO-type procedure \cite{Dumbser:2007a},  in which the reconstructed polynomials are defined for the entire volume 
$\Omega_i$,  as detailed in \cite{castro2012}.  Recall that for first-order Godunov-type methods there are many ways of solving the conventional, piece-wise constant Riemann problem.  The situation for the high-order ADER schemes is analogous; there are by now several ways of solving the associated {\bf Generalised Riemann Problem, or GRP} \cite{Toro:2020a}. In the GRP the data is no longer piece-wise constant, but piece-wise smooth,  given for example by polynomials of arbitrary degree. Moreover, in the GRP, the source terms in the equations are included.
Of the several solvers available for the generalised Riemann problem we adopt the one fully described in \cite{castro_crist:2008}.  More specifically, for the 2D shallow water equations of intereset here we follow \cite{castro2012}.  Details are omitted. 

\section{One-dimensional tests for the first-order scheme}

In this  section we consider one-dimensional test problems to assess the novel, first-order splitting scheme 
proposed in this paper. Three  Riemann problems and one steady problem for flow over a bump are solved.  All four problems have reference solutions, which are used to assess the flux splitting numerical solutions.  For the assessement we also include solutions from the Godunov upwind method used in conjunction with exact Riemann solver \cite{Toro:2001a}.  For the steady-state test we run simulations in both one and two space dimensions.   
Apart from the stationary shock over a bump, all simulation of this section  are first-order accurate, in which the scheme uses numerical fluxes (\ref{eq:1D11}) and  (\ref{eq:1D12}),  in conjunction with the two-rarefaction approximation (\ref{eq:1D29}).  The simulations for the stationary shock over a bump are second-order accurate and well-balanced.

\subsection{Riemann problems}

\begin{figure}
      \centerline{
      \includegraphics[width=0.8\columnwidth]{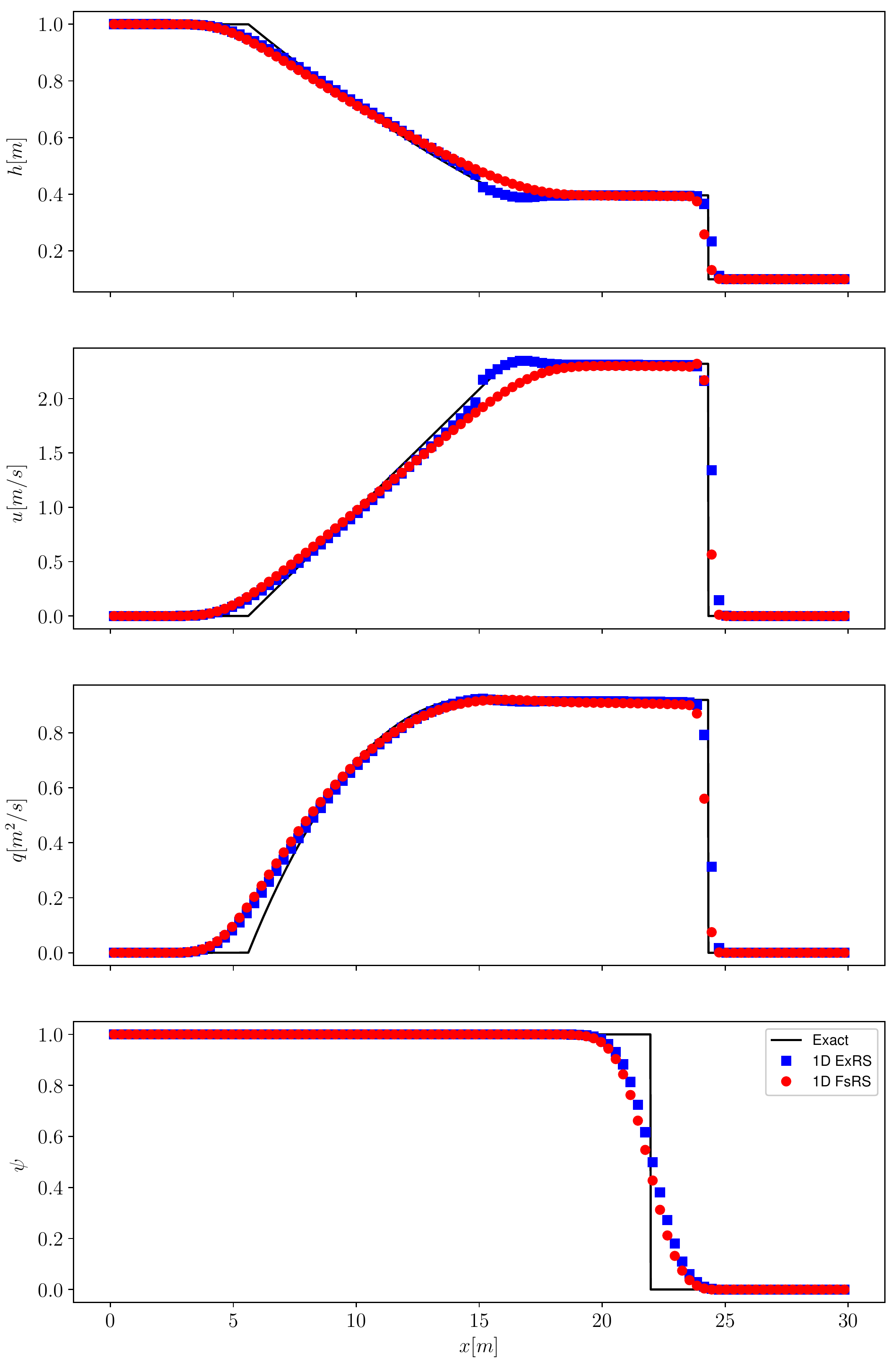}
       }
       \caption{\textbf{Results for Test 2: left shock, middle contact and right shock}.  Comparison of  numerical results from the present TV-type flux vector splitting scheme against the exact solution and the numerical solution from the Godunov upwind method used in conjunction with the exact Riemann solver.
Mesh used $M$ = 100; domain length $L$=30 m, output time $T_{out}$ = 3 s and CFL coefficient $C_{cfl}$ = 0.9.}
\label{fig:Test1}
\end{figure}

\begin{figure}
      \centerline{
      \includegraphics[width=0.8\columnwidth]{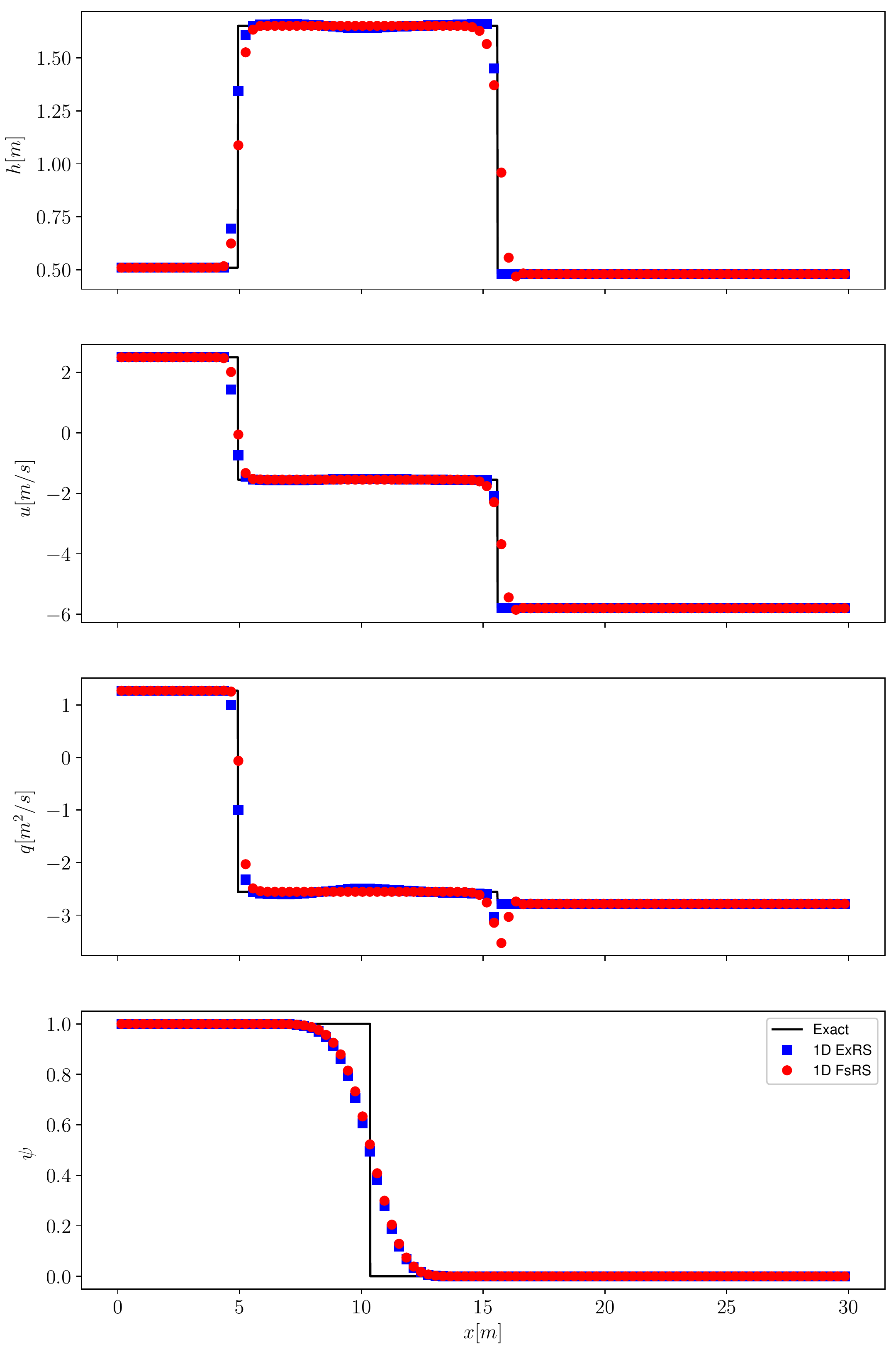}
       }
\caption{\textbf{Results for Test 2: left shock, middle contact and right slowly moving shock}.  Comparison of  numerical results from the present TV-type flux vector splitting scheme against the exact solution and the numerical solution from the Godunov upwind method used in conjunction with the exact Riemann solver.
Mesh used $M$ = 100; domain length $L$=30 m, output time $T_{out}$ = 3 s and CFL coefficient $C_{cfl}$ = 0.9.}
\label{fig:Test2}
\end{figure}

\begin{figure}
      \centerline{
      \includegraphics[width=0.8\columnwidth]{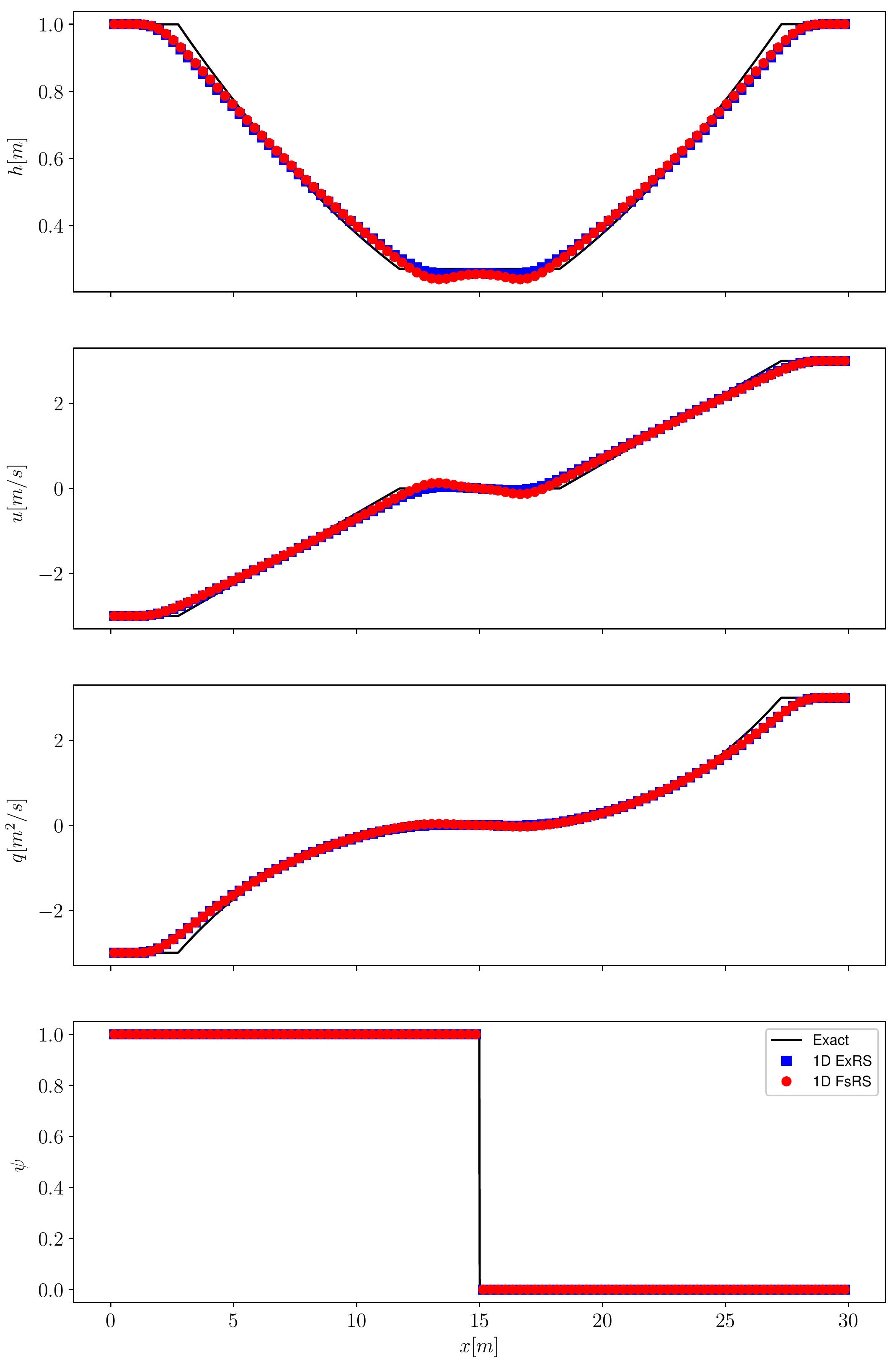}
       }
\caption{\textbf{Results for Test 3: left rarefaction, stationary contact and right rarefaction}.  Comparison of  numerical results from the present TV-type flux vector splitting scheme against the exact solution and the numerical solution from the Godunov upwind method used in conjunction with the exact Riemann solver.
Mesh used $M$ = 100; domain length $L$=30 m, output time $T_{out}$ = 2 s and CFL coefficient $C_{cfl}$ = 0.9.}
\label{fig:Test3}
\end{figure}

The initial conditions for the initial left and right data are given in Table~\ref{tab:initial_RP}. The spatial domain is defined by the interval $[0, 30]$, in meters.  For all the tests, the initial discontinuity is located in the middle of the computational domain, at $x$=\SI{15}{m}.   The mesh used has $M=100$ cells and the CFL coefficient is $C_{cfl}=0.9$.  In each figure {\bf ExRS} means Godunov's method used with the exact solver for the shallow water Riemann problem, while {\bf FsRS} means the present TV-type flux vector splitting scheme.  
Figs.~\ref{fig:Test1} to \ref{fig:Test3}  show solution profiles at a fixed time for water depth, velocity, discharge and concentration of the passive scalar across the complete wave structure.  

Fig.~\ref{fig:Test1} shows the solution for Test 1, which  consists of a left rarefaction, a contact discontinuity and a right shock. The results from our flux-splitting scheme are, overall, comparable with those from Godunov's method used in conjunction with the exact Riemann solver.  The well-known entropy glitch afflicts the Godunov method but not our method, at the cost, however, of resolving the tail of the left rarefaction less accurately.

Fig.~\ref{fig:Test2} shows results for Test 2, the solution of which  consists of three discontinuities: two shocks and one contact wave. The left shock and the contact wave travel to the left; the right shock has a small positive speed. This is a severe test for which many numerical methods fail due to the presence of the right slowly moving shock, which tends to generate visible oscillations, even for first-order methods. 
In comparison with Godunov's method, the shock waves are slightly more diffused.  Behind the right shock the Godunov's method has small oscillations that are not observed in our method.  For the contact
wave the resolution is comparable.  The discharge plot shows a kind of undershoot at the right slowly moving shock, which is more evident for the flux-splitting method of this paper.

Fig. \ref{fig:Test3} shows results for Test 3, whose solution consists of two strong rarefaction waves travelling in opposite directions and one stationary contact discontinuity. The water in the star region between the two rarefaction waves is very shallow, which may cause some schemes to  fail after computing negative depth values. In the results from the flux-splitting scheme, the resolution of the two rarefaction waves is comparable with that of Godunov's method and the stationary contact discontinuity is sharply resolved.  In the middle region, in the vicinity of the stationary contact, the Godunov method with the exact Riemann solver gives a more accurate solution than that from the present scheme.\\

Overall, in spite of its simplicity,  our flux-splitting scheme when applied to Riemann problems gives results that are comparable with those from the more sophisticated Godunov's upwind  scheme used in conjunction with the exact Riemann solver. 

\begin{table}
\caption{Initial conditions for Riemann problem tests.}
   \begin{center}
   \begin{tabular}{c|c|c|c|c|c|c}
   \hline\hline
       test  & $h_L$ [m]   &   $u_L$   [\si{ms^{-1}}] & $\psi_L$ [-] &  $h_R$ [m]   &   $u_R$   [\si{ms^{-1}}] & $\psi_R$ [-] \\ [3pt]    
       \hline    
        1  & 1.0 &  0.0 & 1.0  &  0.1  & 0.0 & 0.0 \\
        2  & 0.51 &  2.5 & 1.0  &  0.48  & -5.8 & 0.0 \\
        3  & 1.0 &  -3.0 & 1.0  &  1.0  & 3.0 & 0.0 \\
        \hline\hline
   \end{tabular}
   \label{tab:initial_RP}
   \end{center}
 \end{table}


\subsection{Transonic flow over a bump with stationary shock}

This problem is designed to test our numerical method for a problem involving steady flow over a bump that contains a stationary shock. The spatial domain is defined by the interval $[0,25]$, with a variable bathymetry given as
\begin{equation}
    b(x) = \left\{  \begin{array}{ccc}
        0.2 - 0.05 (x-10)^2 & \mbox{if} & 8 < x < 12 \,,  \\
        0.0  & \mbox{if} & \mbox{otherwise}\,. 
    \end{array}\right.
\end{equation} 
We compute a steady solution to  this problem, which is determined by the boundary conditions. In this test problem, subcritical flow enters the domain from the left, accelerates up to a critical condition at $x=10$ and reaches a supercritical state until a stationary hydraulic shock is created, returning to a subcritical state downstream. At the steady shock a subcritical regime is recovered. The boundary conditions are defined by subcritical conditions at $x$=\SI{0.0}{m} with $q(x=0,t)$=\SI{0.18}{m^2s^{-1}} and $x$=\SI{25}{m}, with $h(x=25,t)$=\SI{0.33}{m}. The initial conditions used are $h(x,t=0)$=\SI{0.33}{m} and $q(x,t=0)$=\SI{0.18}{m^2s^{-1}}.

\begin{figure}
      \centerline{
      \includegraphics[width=0.9\columnwidth]{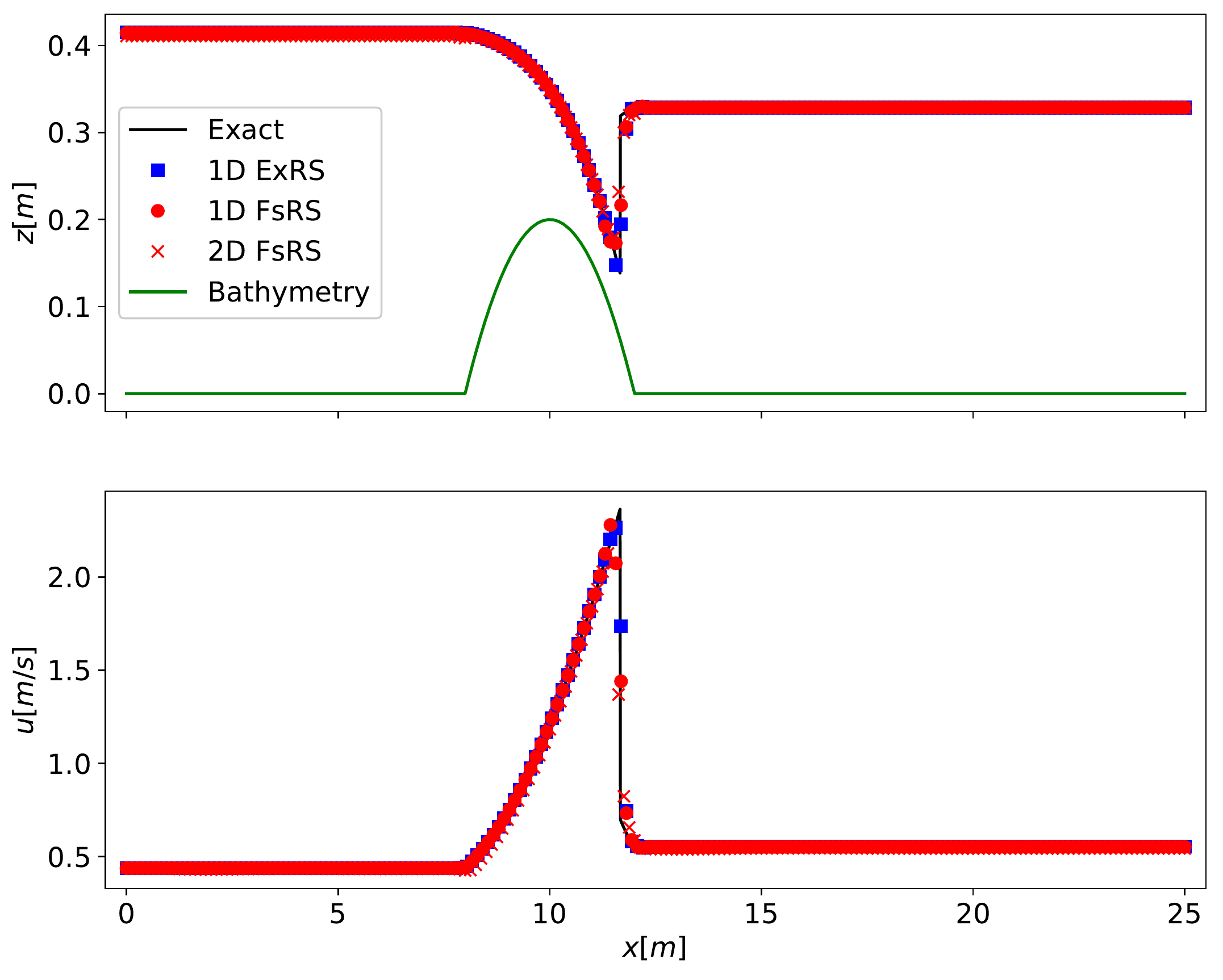}
       }
\caption{\textbf{Transcritical  flow over a bed bump: steady state test}.  Numerical solutions are compared with the reference solution.
Computational parameters are: 1D domain length $L$=25 m, 2D rectangular domain size $L$=25 m $\times$ $W$=25/20 m, simulation time $T_{out}$ = 200 s and $C_{cfl}$ = 0.9.  Number of cells $M_{1D}$ = 200 for the 1D and $M_{2D}$=4420 for 2D (equivalent to an averaged triangle edge size of $L/M_{1D}$=25/200 m),}
\label{fig:StatShock}
\end{figure}

Fig.~\ref{fig:StatShock} shows numerical results from  1D and 2D after a simulation time of $T_{out}$=\SI{200}{s},  time at which the numerical solution has reached steady state conditions.  The black full line depicts the reference solution which is obtained by solving numerically the system of ordinary differential equations resulting from the steady-state equations. The bathymetry depicting the bump is shown by the full green line.  The blue squared symbols identify the results from Godunov's method,  used in conjunction with the exact Riemann solver;  the red full circles show results from the present TV-type flux vector splitting scheme.  The red cross shows results from the 2D splitting scheme on an unstructured triangular mesh.
As for the Riemann problem tests, results for this problem show that our flux-splitting scheme compares very favorably with the reference solution and that from the Godunov method used in conjunction with the exact Riemann solver.\\


%
%
%

\section{Two dimensional tests for high-order schemes}

In this section we test the high-order extension of the present TV-type flux vector splitting scheme. To this end we use the high-order ADER finite volume framework \cite{Toro:2001}.  We use the HEOC Generalized Riemann solver from \cite{castro_crist:2008} and the well-balanced approach presented in \cite{castro2012}.  Details are omitted, as they have been published in the references above.  First we verify the method's accuracy, then we simulate a circular dam break problem and finally we simulate the propagation of a tsunami wave in a realistic scenario over the Pacific Ocean. 

\subsection{Convergence rates  study}

\begin{table}
\begin{center}
\begin{tabular}{ c | c c | c c | c c}
\hline
$\Delta x$ & Error & $L_1$ & Error & $L_2$ & Error & $L_\infty$ \\
\hline
0.109 & 7.51e-03 &  2.36   &    5.36e-03 &  2.83   &    1.92e-02 &  4.34 \\
0.080 & 3.37e-03 &  2.51   &    2.78e-03 &  2.42   &    1.98e-02 &  1.69 \\
0.069 & 1.14e-03 &  3.56   &    8.35e-04 &  3.68   &    5.06e-03 &  3.35 \\
0.059 & 2.71e-04 &  4.68   &    2.36e-04 &  4.55   &    1.04e-03 &  4.68 \\
0.048 & 9.41e-05 &  4.79   &    7.77e-05 &  4.75   &    3.68e-04 &  4.78 \\
0.046 & 7.08e-05 &  4.82   &    6.72e-05 &  4.65   &    2.82e-04 &  4.79 \\
0.034 & 1.58e-05 &  4.90   &    1.27e-05 &  4.88   &    5.24e-05 &  5.01 \\
0.024 & 2.23e-06 &  5.07   &    1.93e-06 &  5.01   &    9.65e-06 &  5.00 \\
\hline
\end{tabular}
\caption{\textbf{Convergence-rate study in 2D}.  ADER scheme in conjunction with the TV-type flux vector splitting scheme of this paper.  Results for the 5\textsuperscript{th} order scheme are shown for the output time $T_{out}$=\SI{1}{s}. Convergence rates are calculated for the discharge per unit width $q_x$.  The fifth order of accuracy is attained in all norms for the finest mesh.
}\label{tab:ConvRate}
\end{center}
\end{table}

Here we verify that the expected high accuracy,  up to fifth order in space and time, of our numerical scheme  is actually attained; we do so by studying empirical convergence rates.  For the assessment we compare the numerical solutions against a {\it manufactured}, prescribed,  exact solution presented in \cite{castro2012}. The error is measured in three different norms. The manufactured exact solution $\mathcal{Q}(\vec{x},t)$ is given by
\begin{equation}\label{eq:ConvExacSol}
    \mathcal{Q}(\vec{x},t) = \left[ \begin{array}{c}
    h(\vec{x},t) \\
    h(\vec{x},t) u(\vec{x},t) \\
    h(\vec{x},t) v(\vec{x},t)
    \end{array} \right]\mbox{where} \begin{array}{rcl}
    h(\vec{x},t) & = & e^{0.1t} - b(\vec{x}) \\
    u(\vec{x},t) & = & 0.2+0.1\sin(x\pi) \\
    v(\vec{x},t) & = & 0.2+0.1\sin(y\pi) \\
    b(\vec{x}) & = & 0.2\,e^{-8(x^2+y^2)}  \;,
    \end{array}
\end{equation}
with $\vec{x}=(x,y)$. Of course, the proposed function (\ref{eq:ConvExacSol}) is not an exact solution of the 
original equations (\ref{eq:2D_vectorial}). Therefore, when $\mathcal{Q}(\vec{x},t)$ is plugged into  (\ref{eq:2D_vectorial}), a new source term appears that restores the balance, namely  $\mathcal{S}(\mathcal{Q}(\vec{x},t)) = \partial_{t}{\mathcal{Q}} + \partial_{x}{\bf F}_x({\mathcal{Q}}) + \partial_{y}{\bf F}_y({\mathcal{Q}}) - {\bf S}$.  For the convergence rates study we solve the new problem
\begin{equation}\label{eq:ConvTest}
\begin{array}{ll}
    \mbox{PDE:} & \partial_{t}{\bf Q} + \partial_{x}{\bf F}({\bf Q}) + \partial_{y}{\bf G}({\bf Q}) = {\bf S} + \mathcal{S}(\mathcal{Q}(\vec{x},t)) \;, \\
    \mbox{IC:} & {\bf Q} = \mathcal{Q}(\vec{x},t)) \;.
\end{array}
\end{equation}
This means that the tested methods will have to solve 2D problems with source terms to the desired high-order of accuracy. Note that even if the original problem is homogeneous (no source terms) the {\it manufactured solution} satisfies  equations with source terms.
The initial-boundary value problem (\ref{eq:ConvTest}) is defined in a two-dimensional square domain $\Omega \in [0,1] \times [0,1]$,  discretized with a triangular  grid  and solved until time $T_{out}$=1 s. The numerical solution is obtained for a sequence of successively refined meshes.  Then by comparing the numerical error between two successive mesh refinements it is possible to compute the order of convergence of the numerical method.  In Table \ref{tab:ConvRate} the results for the fifth order method are presented, where $\Delta x$ is a measure of the element mesh size. The errors are computed for the specific discharge  $q_x$ in the $x$-direction.   As seen, our method attains the expected convergence rate of fifth order (in space and time) in all norms. Similar results are obtained for second, third and fourth orders.

\begin{figure}
      \centerline{
      \includegraphics[width=0.8\columnwidth]{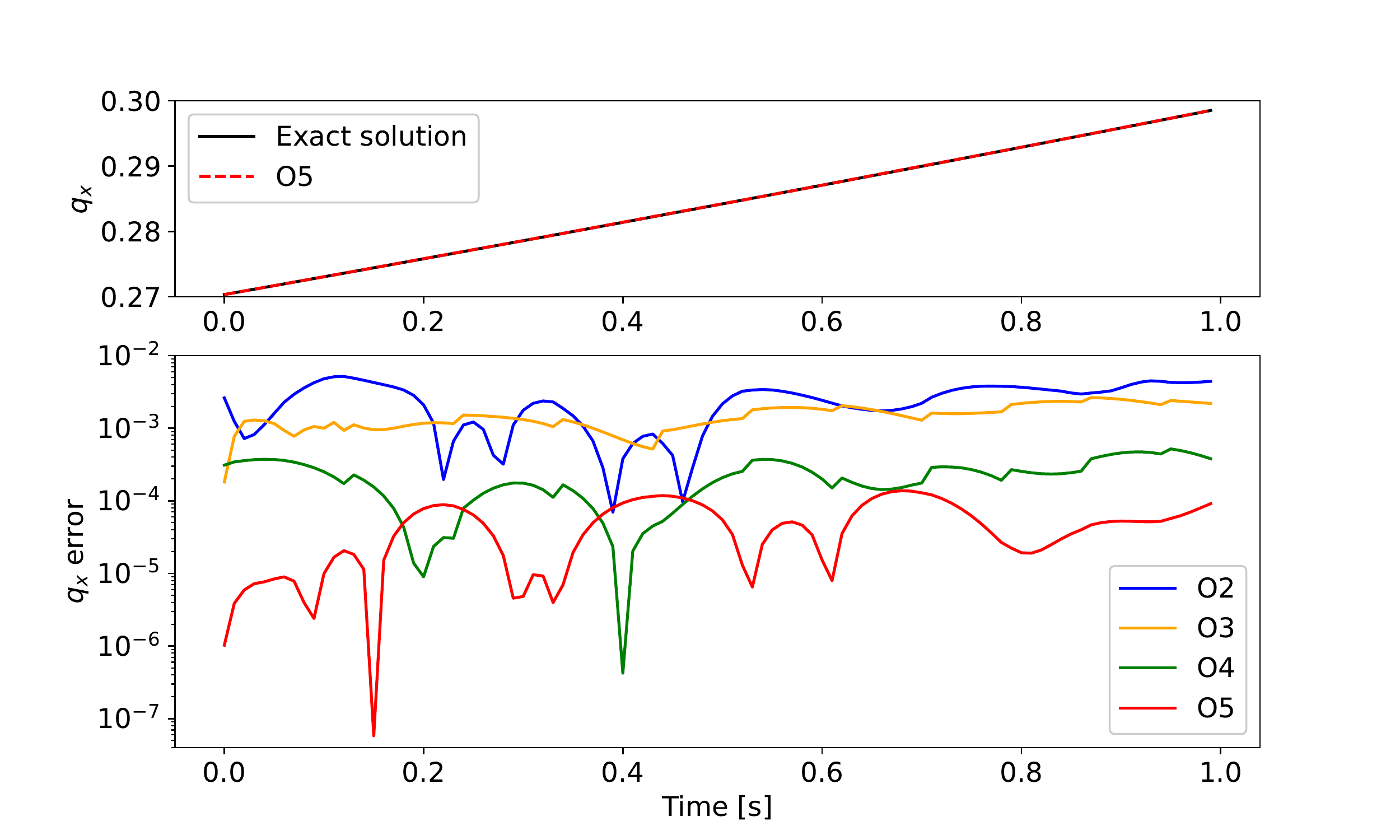}
       }
\caption{\textbf {Time evolution for specific discharge $q_x$ at location $x$=0.25m,$y$=0.75m}.  Top panel shows the exact solution from  (\ref{eq:ConvExacSol}) (black line) compared to the fifth order accurate numerical solution (dashed blue line).  The bottom panel shows the local error of the solutions obtained from the second to fifth order schemes. }
\label{fig:ConvergenceTest_Pickpoints}
\end{figure}

Fig.~\ref{fig:ConvergenceTest_Pickpoints} depicts the time evolution of the numerical solution for  $q_x(x,y,t)$ at the spatial position $x$=0.25m, $y$=0.75m, from the initial time $t=0$ s until time $t=1.0$ s.  In the upper panel we show the exact solution from  (\ref{eq:ConvExacSol}) (black line) and the fifth order accuracy numerical solution (dashed blue line). In the bottom panel we show the error defined as $err = |hu - hu_{ex}|/hu_{ex}$ for the solutions obtained from the second to fifth order schemes.  We observe how the error decreases as the order of accuracy of the method increases, even if we see specific time points at which low-order methods exhibit smaller local errors. The dominant behavior is the one expected, where the high-order methods used with the new splitting method strongly reduces the error.

\subsection{Circular-dam break test problem}

\begin{figure}
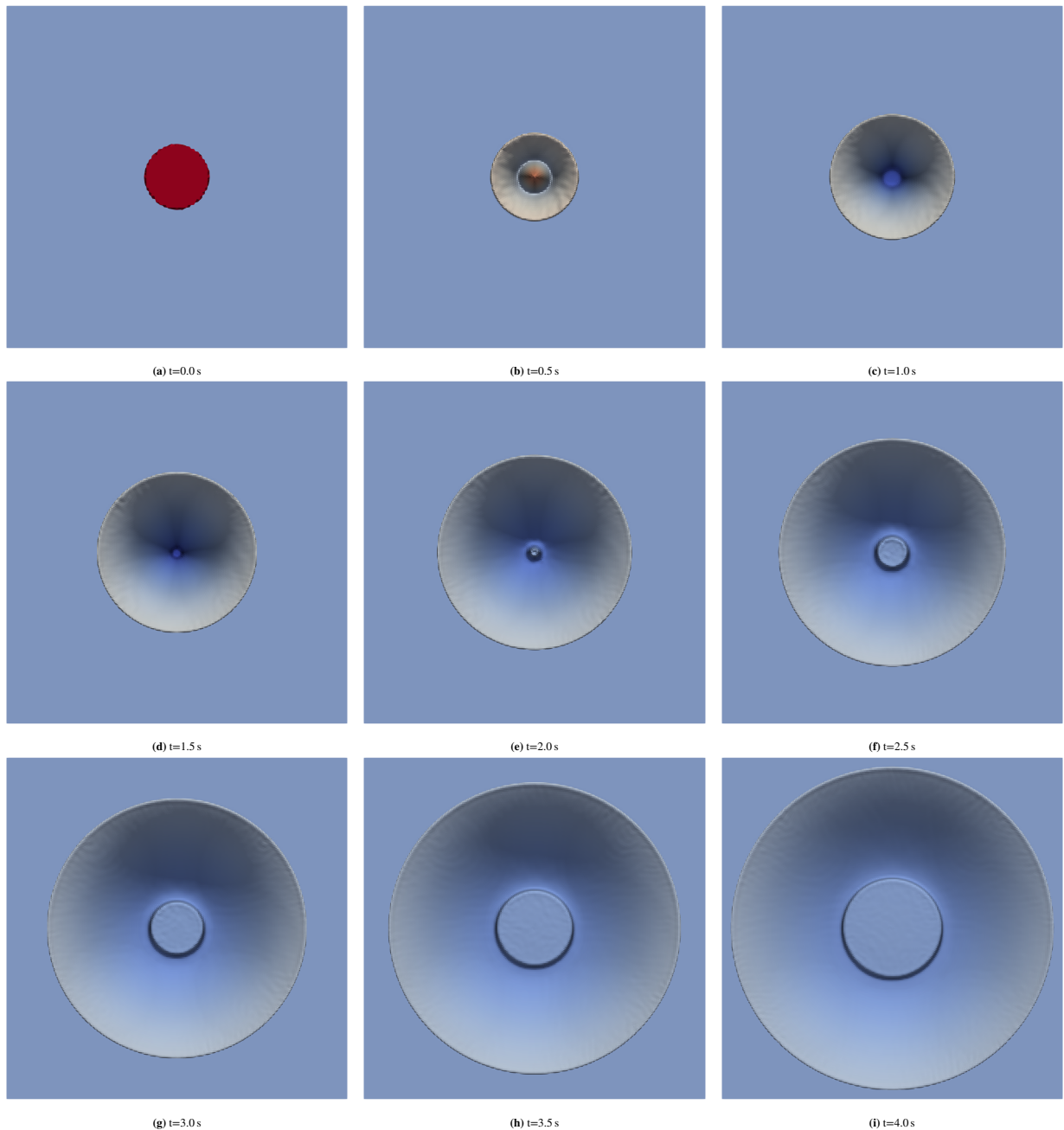

      \centering
      \begin{subfigure}[b]{0.3\textwidth}
      \includegraphics[width=\textwidth]{figs/SnapShot_0.0s}
      \caption{t=\SI{0.0}{s}}
      \end{subfigure}
      \begin{subfigure}[b]{0.3\textwidth}
      \includegraphics[width=\textwidth]{figs/SnapShot_0.5s}
      \caption{t=\SI{0.5}{s}}
      \end{subfigure}
      \begin{subfigure}[b]{0.3\textwidth}
      \includegraphics[width=\textwidth]{figs/SnapShot_1.0s}
      \caption{t=\SI{1.0}{s}}
      \end{subfigure}
      \begin{subfigure}[b]{0.3\textwidth}
      \includegraphics[width=\textwidth]{figs/SnapShot_1.5s}
      \caption{t=\SI{1.5}{s}}
      \end{subfigure}
      \begin{subfigure}[b]{0.3\textwidth}
      \includegraphics[width=\textwidth]{figs/SnapShot_2.0s}
      \caption{t=\SI{2.0}{s}}
      \end{subfigure}
      \begin{subfigure}[b]{0.3\textwidth}
      \includegraphics[width=\textwidth]{figs/SnapShot_2.5s}
      \caption{t=\SI{2.5}{s}}
      \end{subfigure}
      \begin{subfigure}[b]{0.3\textwidth}
      \includegraphics[width=\textwidth]{figs/SnapShot_3.0s}
      \caption{t=\SI{3.0}{s}}
      \end{subfigure}
      \begin{subfigure}[b]{0.3\textwidth}
      \includegraphics[width=\textwidth]{figs/SnapShot_3.5s}
      \caption{t=\SI{3.5}{s}}
      \end{subfigure}
      \begin{subfigure}[b]{0.3\textwidth}
      \includegraphics[width=\textwidth]{figs/SnapShot_4.0s}
      \caption{t=\SI{4.0}{s}}
      \end{subfigure}
\caption{\textbf{Time evolution  of the circular dam break problem}. Solutions for depth are shown from time $t=0.0$ to $t=4.0$ every $0.5$ seconds.  Red shows deeper water, while the dark blue represents shallow water zones.}
\label{fig:CircularDam_snapshots}
\end{figure}

In this section we present numerical results for a 2D circular dam break problem, for which an accurate, radial quasi 1D numerical solution can be computed and used as a reference solution.  See Chap.  13 of \cite{Toro:2001a} for details on circular dam break problems.  For our test problem here the physical domain is defined by the square $0 \le x \le L$, $0 \le  y \le L$, with $L$=\SI{40}{m}.  The bottom distribution is $b(x,y)=0$ throughout.The initial conditions at time $t=0$ are defined by two constant states separated by the circle $r(x,y)=\sqrt{(x-20)^2+(y-20)^2}  = 2.5$.  In the inner circular region the  water depth is higher than that in the outer circular region,  that is $h_{in}(x,y,0)=2.5 m$ and $h_{out}(x,y,0)=1.0 m$.  The velocity at the initial time is zero everywhere. The  initial condition for depth is graphically shown in the top left corner of  Fig.~\ref{fig:CircularDam_snapshots} in which the red color represents the higher water depth $h=2.5$ and the light blue colour represents $h=1.0$.  The physical domain is discretized by a mesh of $92258$ triangles. \\

 Fig.~\ref{fig:CircularDam_snapshots}  shows water depth $h(x,y, t)$ at fixed times, starting from the initial condition at $ t=0$, until time $t=4.0$, every 0.5 seconds.  The instantaneous  collapse of the circular dam generates two main waves, namely an outward travelling circular shock wave called the primary shock, and an inner travelling circular rarefaction wave.  The primary shock travels outwardly with decreasing speed and shock strength. The inner facing rarefaction wave reaches the centre and reflects as a rarefaction wave,  strongly reducing the water depth. The reflected rarefaction overexpands the flow causing the formation of a secondary,  inner-facing shock wave. This secondary shock is initially swept outwardly for a short time, then comes to a halt, to then begin to propagate towards the centre.  The inner shock implodes into the centre and reflects as a shock wave,  obviously,  and travels outwardly  behind the primary shock.  For a full discussion of the evolving wave patterns see Chap.  13 of \cite{Toro:2001a}.

\begin{figure}
      \centerline{
      \includegraphics[width=0.5\columnwidth]{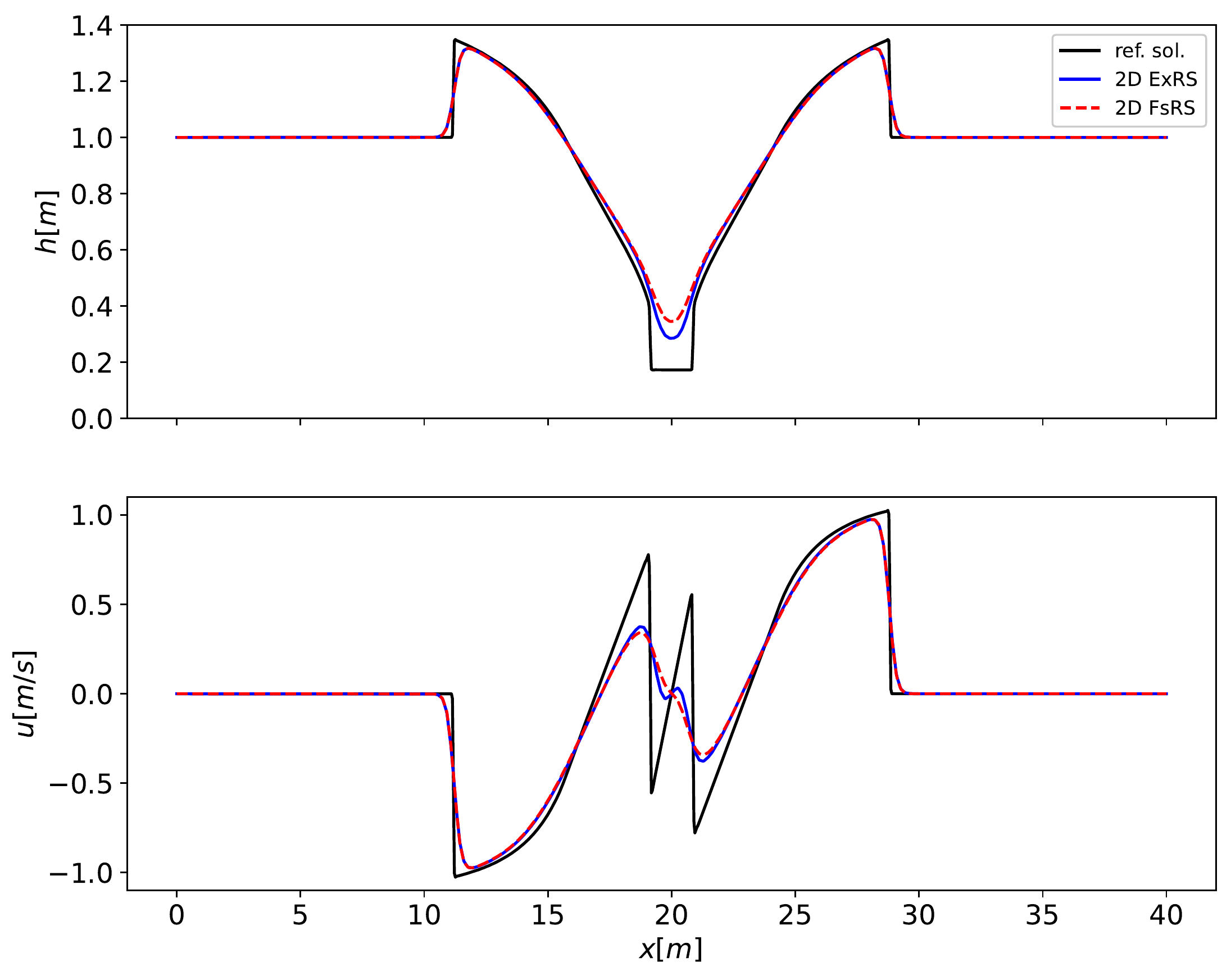}
      \includegraphics[width=0.5\columnwidth]{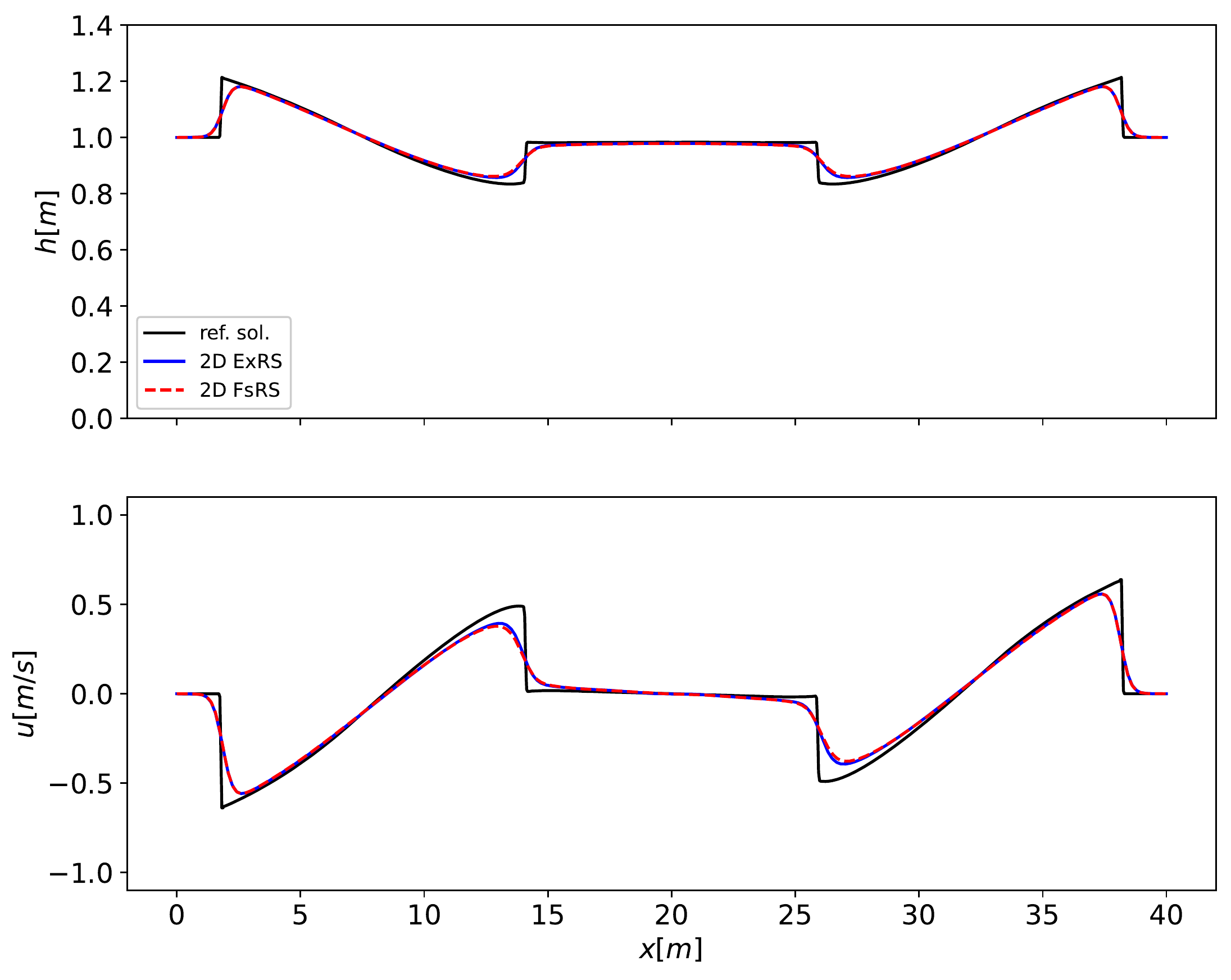}
       }
\caption{\textbf{Circular dam break problem: first-order numerical results}.  Solution is displayed at  time 1.5s (left panels) and 4.0s (right panels) for water depth (upper panels) and particle velocity (lower panels).}
\label{fig:CircularDam_p0}
\end{figure}

\begin{figure}
      \centerline{
      \includegraphics[width=0.5\columnwidth]{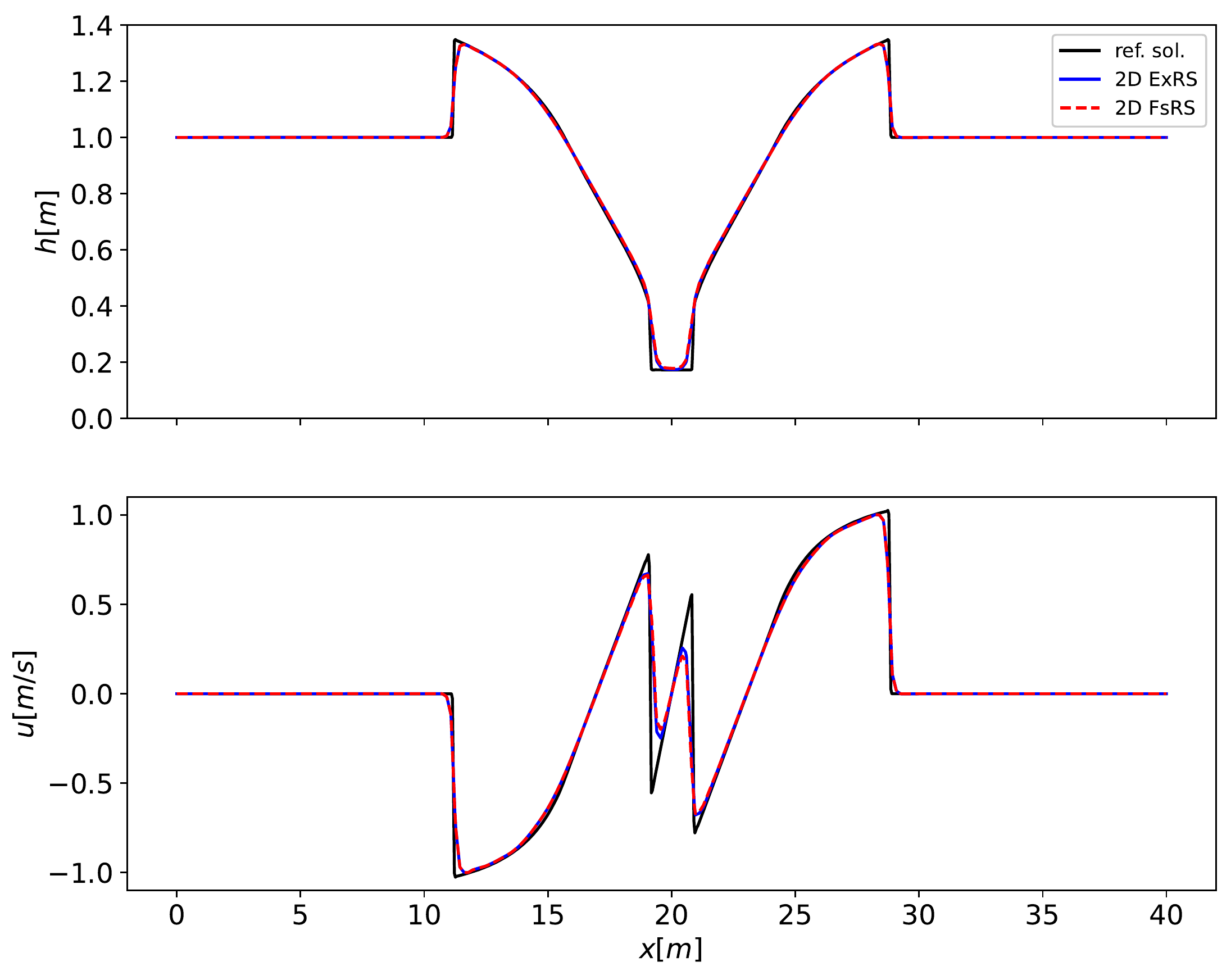}
      \includegraphics[width=0.5\columnwidth]{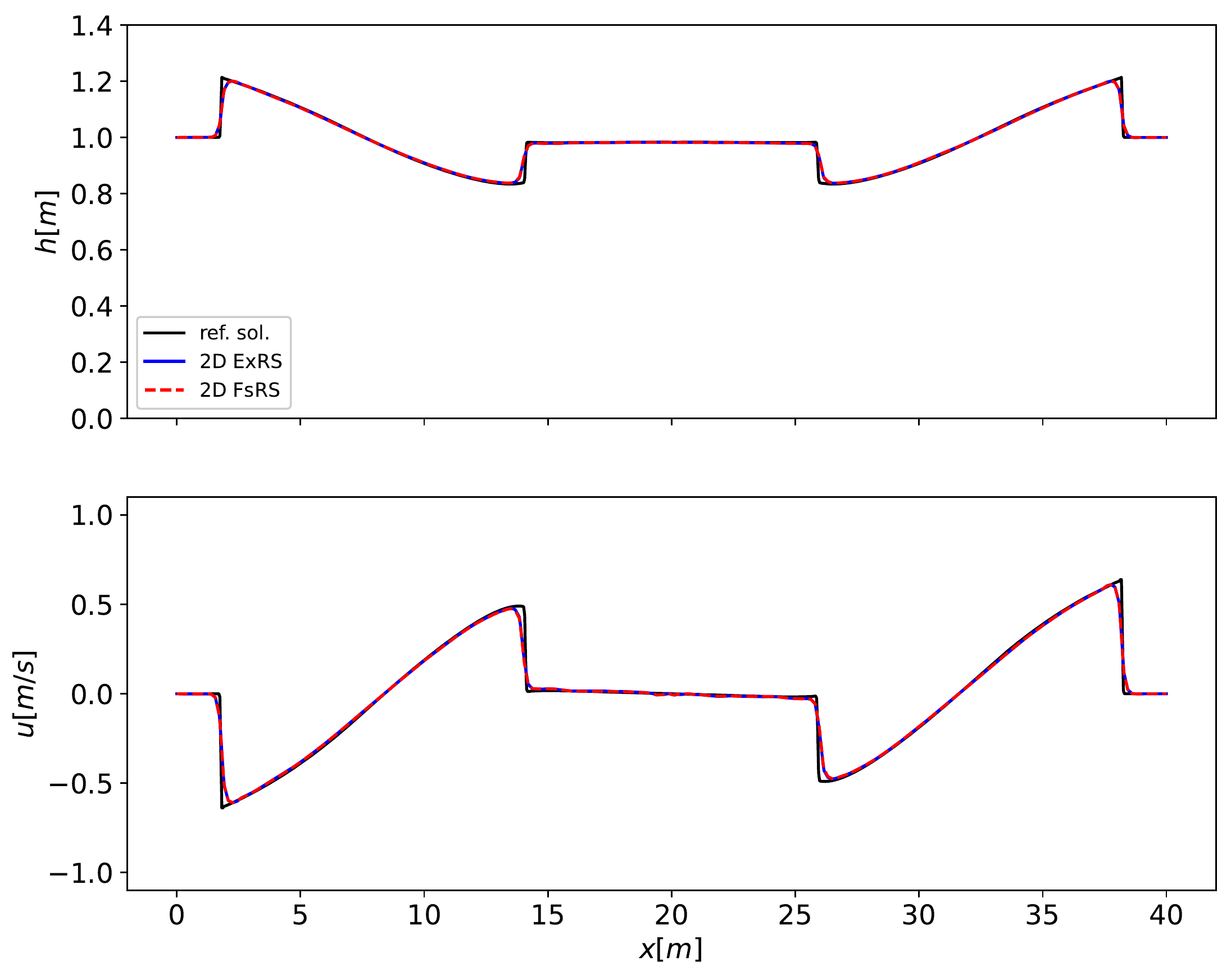}
       }
\caption{\textbf{Circular dam break problem: second order numerical results}. Solution is displayed at  time 1.5s (left panels) and 4.0s (right panels) for water depth (upper panels) and particle velocity (lower panels).}
\label{fig:CircularDam_p1}
\end{figure}

\begin{figure}
      \centerline{
      \includegraphics[width=0.5\columnwidth]{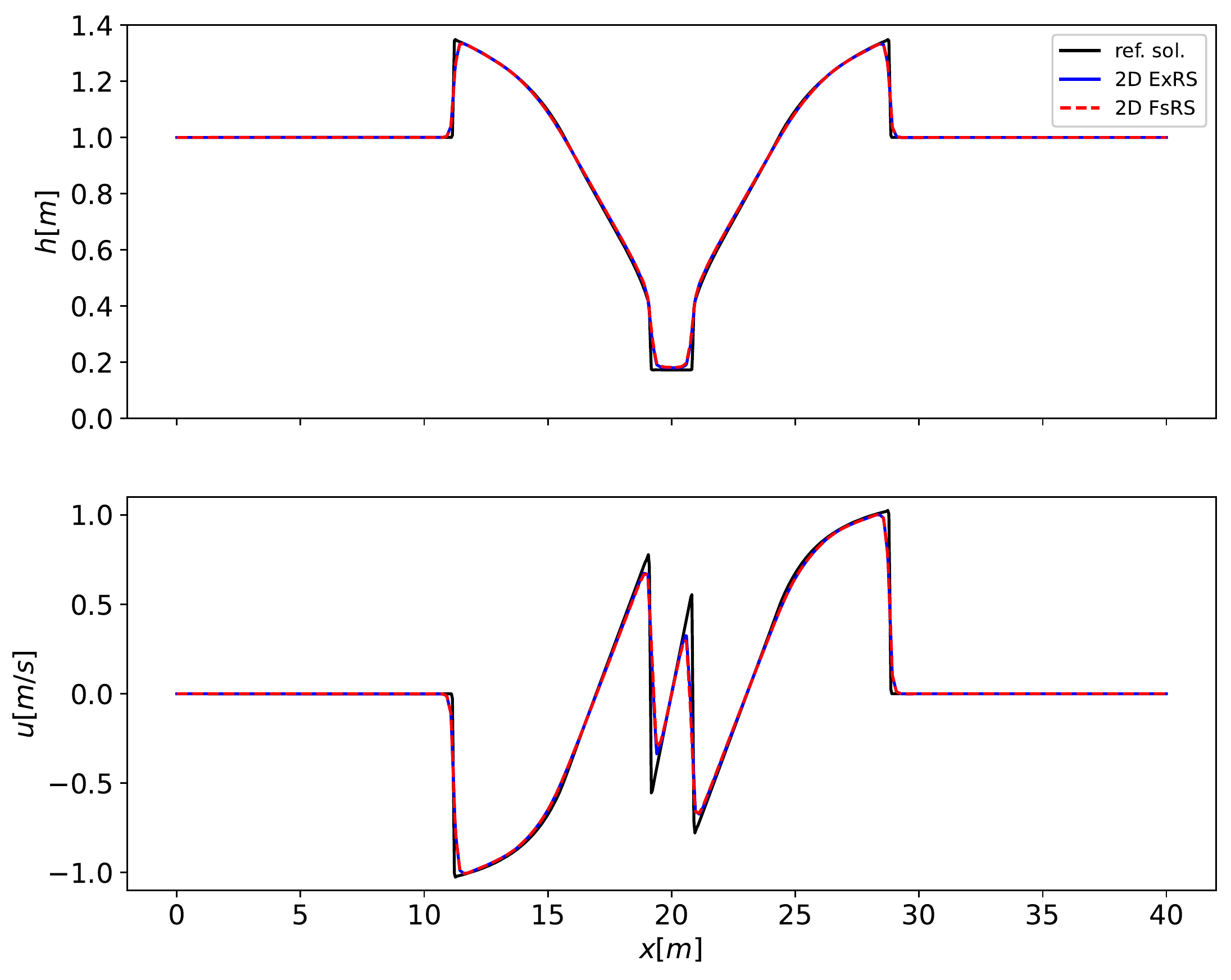}
      \includegraphics[width=0.5\columnwidth]{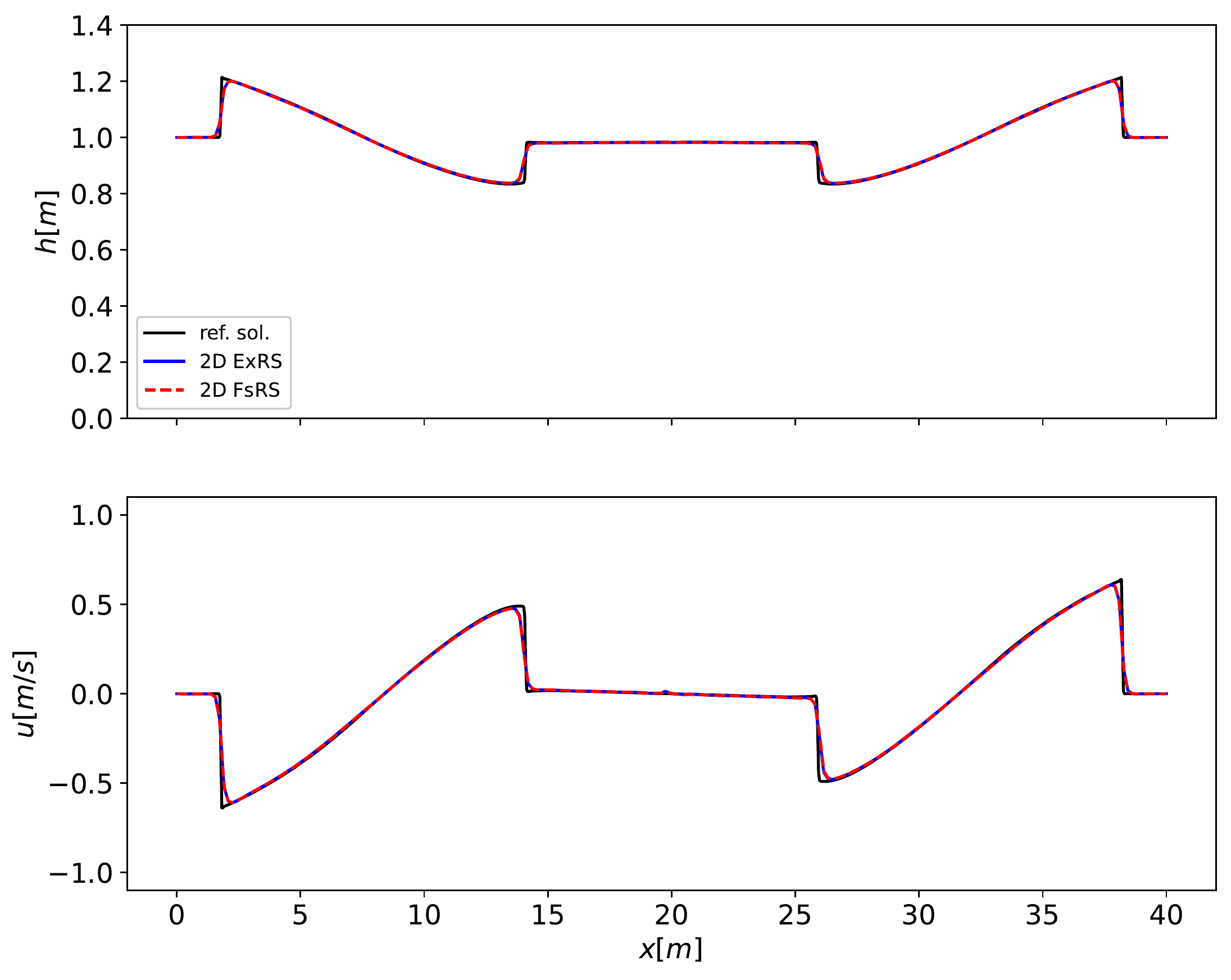}
       }
\caption{\textbf{Circular dam break problem: third order numerical results}. Solution is displayed at  time 1.5s (left panels) and 4.0s (right panels) for water depth (upper panels) and particle velocity (lower panels).}
\label{fig:CircularDam_p2}
\end{figure}

Figs.~\ref{fig:CircularDam_p0}  to ~\ref{fig:CircularDam_p2} present numerical results for a slice cut along the line $0<x<40$ and $y=20$. Results are shown for the first to third order schemes; the solid blue line shows the solution from the Godunov method in conjunction with the exact Riemann solver, while the dashed red line shows that from the TV-type splitting scheme of this paper. These numerical solutions are compared with a reference solution shown by solid black line,  computed with a 1D radial symmetry finite volume code from NUMERICA Library \cite{Toro:1999c}. More details are found in Chap.  13 of \cite{Toro:2001a}.

We observe that the primary outward-facing shock is resolved with the same accuracy by both numerical methods, see Fig.~\ref{fig:CircularDam_p0}. The same behaviour is reproduced across the rarefaction wave. 
By examining the particle velocity at time $t=1.5$ it is evident that the splitting solver is more diffusive than the Godunov upwind method with the exact Riemann solver.  At the later time $t=4.0$ this difference has become smaller.  For second order Fig.~\ref{fig:CircularDam_p1} and third order Fig.~\ref{fig:CircularDam_p2} the differences between the two numerical solutions are almost negligible; see for example the particle velocity at time $t$=1.5s and $t$=4.0s in Fig.~\ref{fig:CircularDam_p2}. 

\subsection{Tsunami Wave Propagation}

The purpose of this test problem is to assess the practical capability of the numerical methods presented in this paper,  to solve problems of interest to the environmental community. To this end we consider the propagation of a tsunami wave,  in a realistic scenario, over the Pacific Ocean, in the coastal region near the city of Arica,  Chile. The bathymetry is obtained from the General Bathymetric Chart of Oceans (GEBCO) \cite{Gebco:2009} database and the initial condition is an assumed Gaussian function defined at the left boundary of the computational domain to simulate a tsunami wave travelling to the right, towards the western coast of the South American continent.  

\begin{figure}
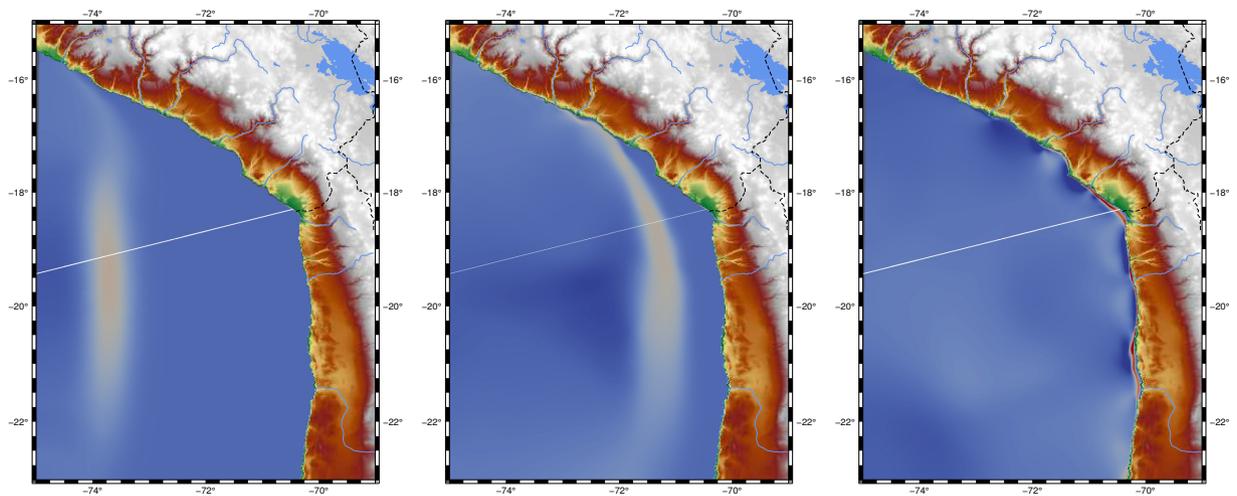

      \centerline{
      \includegraphics[width=0.3\columnwidth]{figs/Tsunami_10min-eps-converted-to}
      \includegraphics[width=0.3\columnwidth]{figs/Tsunami_30min-eps-converted-to}
      \includegraphics[width=0.3\columnwidth]{figs/Tsunami_50min-eps-converted-to}
       }
\caption{\textbf{Tsunami wave propagation test problem}. Wave propagation after 10, 30 and 50 minutes}
\label{fig:TsunamiTest_WavePropagation}
\end{figure}

\begin{figure}
      \centerline{
      \includegraphics[width=0.8\columnwidth]{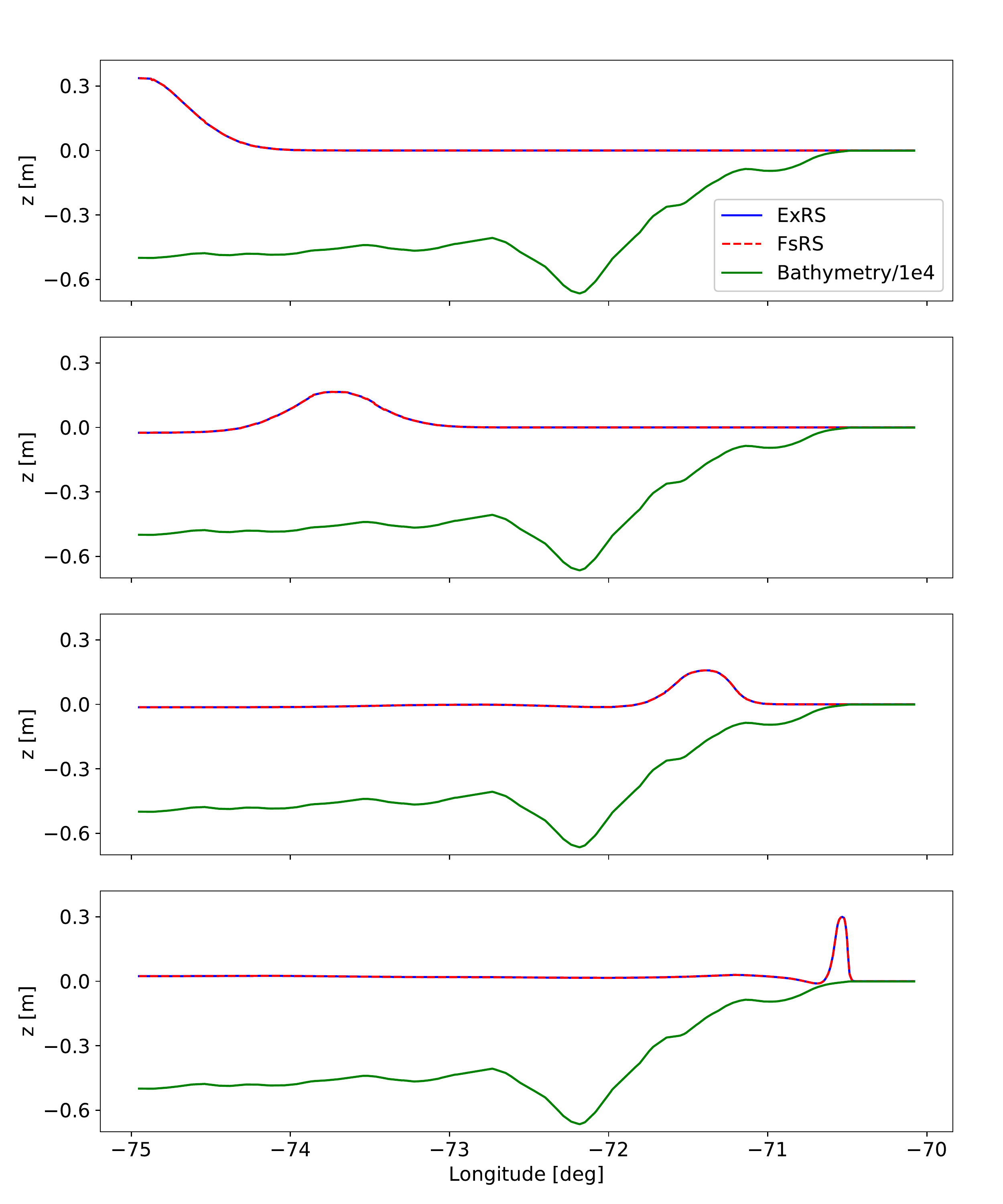}
       }
\caption{\textbf {Tsunami wave propagation test problem: second order numerical solution}. Free surface at 0, 10, 30 and 50 minutes from initial time.}
\label{fig:TsunamiTest_FreeSurface}
\end{figure}

\begin{figure}
      \centerline{
      \includegraphics[width=0.8\columnwidth]{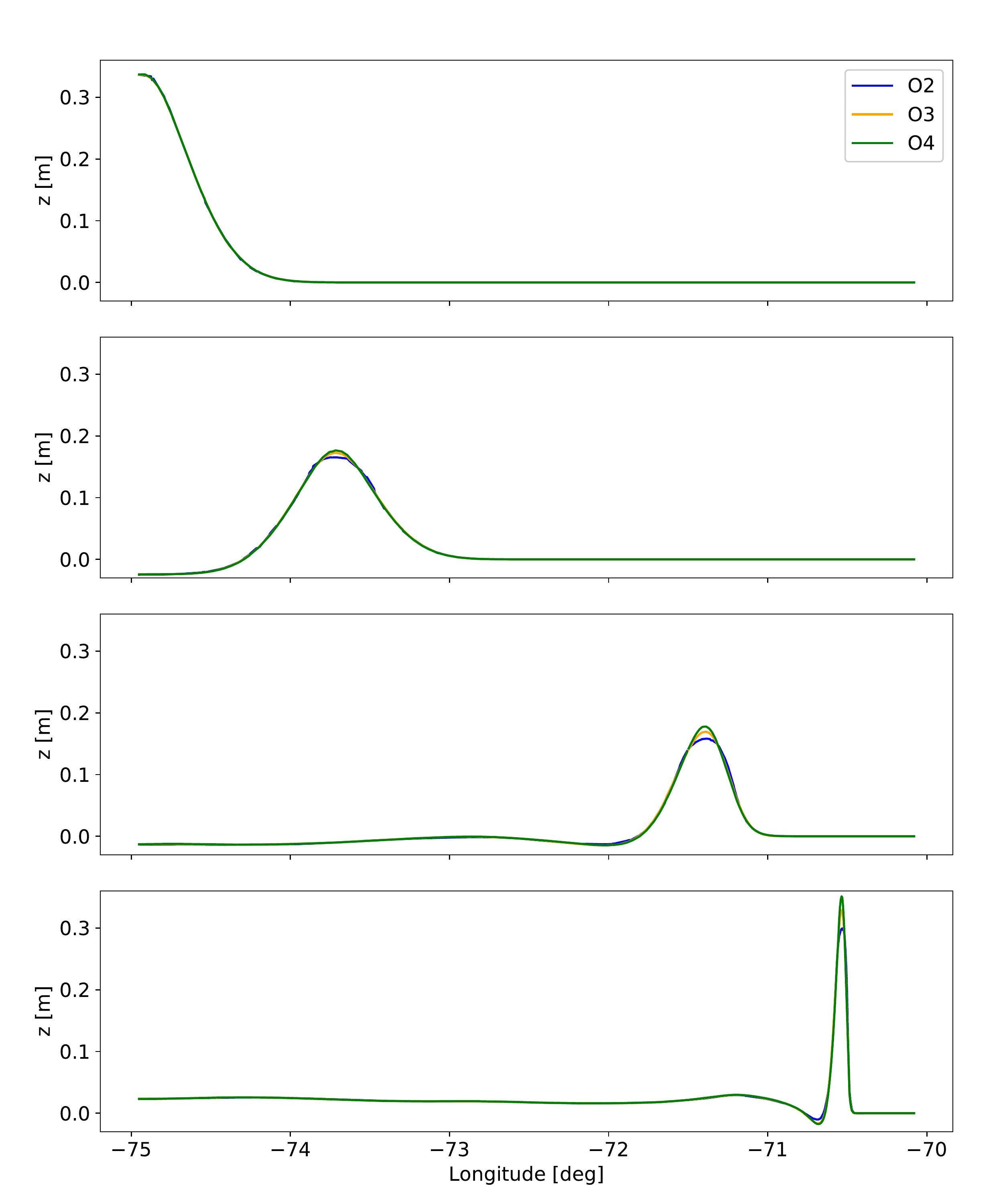}
       }
\caption{\textbf {Tsunami wave propagation test problem: comparison among different high order solutions}. Order 2 (blue line O2), order 2 (yellow line O3), order 4 (green line O4).    Free surface is shown at 0, 10, 30 and 50 minutes from initial time.}
\label{fig:TsunamiTest_FreeSurface_p1p2p3}
\end{figure}

In Fig.~\ref{fig:TsunamiTest_WavePropagation} the tsunami wave propagation is presented after 10, 30 and 50 minutes from the initial time.    We analyze the results in  more details by slicing the domain along the white line in Fig.~\ref{fig:TsunamiTest_WavePropagation} and present the results in terms of free surface elevation in Figs.~\ref{fig:TsunamiTest_FreeSurface} and  \ref{fig:TsunamiTest_FreeSurface_p1p2p3}.  
Note that the bathymetry is plotted using a scaling factor of $10^{-4}$ in order to appreciate the wave amplitude, which is in the order of 20 cm.  In Fig.~\ref{fig:TsunamiTest_FreeSurface_p1p2p3} we compare the numerical solutions from schemes of different orders of accuracy at times $t=10$,  $t=30$  and $t=50$ minutes from initial time. 
It is seen as both the Godunov scheme with the exact Riemann  solver and the TV-type splitting of this paper  agree quite well  on the arrival time and tsunami wave amplitude for the entire time period considered.  Not surprisingly,  the wave amplitude is better preserved by the higher-order methods. This is a crucial observation of relevance when computing solutions for very long time/long distance, where the numerical diffusion of low-order methods may have a dramatic effect on the accuracy of the results at the end of the propagation phase and the initiation of the inundation phase of the tsunami.

\section{Summary and concluding remarks}

We have presented an advection-pressure flux-vector splitting method for the one and two-dimensional shallow water equations following the  TV approach of Toro and V\'azquez  \cite{Toro:2012b}.  The TV splitting  technique splits the full system into two subsystems,  namely an advection system and a  pressure system.  In applying the TV approach to the shallow water equations,  a  modification has been found to be necessary, namely the full inclusion of the continuity equation in the pressure system.  The resulting first-order schemes turn out to be exceedingly simple, with accuracy and robustness comparable to that of the sophisticated Godunov method used in conjunction with the exact Riemann solver. The basic methodology has been extended to 2D geometries using unstructured meshes,  thus forming the building block for the construction of numerical schemes of very high order of accuracy following the ADER approach.   Schemes of up to fifth order of accuracy in space and time have been implemented and tested. The presented numerical schemes have been systematically assessed on a carefully selected suite of test problems, with reference solutions,  in one and two space dimensions.The applicability of the schemes was illustrated through  simulations of tsunami wave propagation in the Pacific Ocean.
Potential extensions of the methodology include the proper treatment of source terms,  specially those of the geometric type,  by making full use of the two sub-systems resulting from the splitting approach.  An attractive future extension would be the devising of semi-implicit schemes so as to exploit the disparity of wave speeds emerging separately from the advection system and the pressure system, as done, for example  by Dumbser and Casulli \cite{Dumbser:2016d}  for the compressible Navier-Stokes equations.  It is also expected that   the splitting approach may prove beneficial in treating viscous terms in the shallow water equations as well as for solving extended shallow water models for applications in oceanography, in which many extra equations may be added to the basic system, with complicated source terms.

\clearpage

\subsection{Bibliography}
\bibliography{SWEbib}






\end{document}